\documentclass[12pt]{iopart}
\def \d{ {\rm d} } \def \e{ {\rm e} }  \def \i{ {\rm i} }

\def\Proof{\noindent\bf{Proof.}\rm\quad}

\newtheorem{theorem}{Theorem}[section]
\newtheorem{lemma}{Lemma}[section]

\newtheorem{remark}{Remark}[section]
\let\oldsection\section
\renewcommand\section{\setcounter{equation}{0}\oldsection}
\renewcommand\theequation{\thesection.\arabic{equation}}
\newcommand\qedhere{\hspace{\fill} $\Box$}

\usepackage{amssymb, bm, iopams, setstack, graphicx, epsfig, subfig}
\usepackage{booktabs} % for much better looking tables
\usepackage{epstopdf}

\begin{document}

%\title[ {\rm D  Zhang } et al]{Direct sampling method for solving an inverse source
%problem for the Helmholtz equation}

\title{Locating multiple multipolar acoustic sources using the direct sampling method}

\author{Deyue Zhang$^1$,  Yukun Guo$^2$, Jingzhi Li$^{3}$ and Hongyu Liu$^4$}

\address{$^1$ School of Mathematics, Jilin University, Changchun, P. R. China \\
$^2$ Department of Mathematics, Harbin Institute of Technology, Harbin, P. R. China\\
$^3$ Department of Mathematics,  Southern University of Science and
Technology, Shenzhen, P.~R.~China\\
$^4$ Department of Mathematics, Hong Kong Baptist University,Kowloon Tong, Hong Kong SAR}
\eads{\mailto{dyzhang@jlu.edu.cn}, \mailto{ykguo@hit.edu.cn}, \mailto{li.jz@sustc.edu.cn}, \mailto{hongyu.liuip@gmail.com}}

\begin{abstract}
This work is concerned with the inverse source problem of locating multiple
multipolar sources from boundary measurements for the Helmholtz equation.
We develop simple and effective sampling schemes for location acquisition of the sources with a single wavenumber. Our algorithms are based on some novel indicator functions whose indicating behaviors could be used to locate
multiple multipolar sources. The inversion schemes are totally ``\emph{direct}'' in the sense that only simple integral calculations are involved in evaluating the indicator functions. Rigorous mathematical justifications are provided and extensive numerical examples are presented to demonstrate the effectiveness, robustness and efficiency of the proposed methods.
\end{abstract}


\section{Introduction}
The inverse source problem concerns with the reconstruction of unknown sources from measured scattering data away from the sources.
It arises in many scientific fields and engineering applications, such as antenna
synthesis \cite{Angel,Devaney,Ramm}, active acoustic tomography \cite{Anastasio,Stefanov}, medical
imaging \cite{Albanese,Ammari,Arridge,Fokas} and pollution source tracing \cite{El,Isakov}.

Over recent years, intensive attention \cite{Angel, Bao2015a, Bao2015b, Devaney,El1,El2,Eller,K1,K2,K3,Ramm}
has been focused on the inverse source problem of determining a source
$F$ in the Helmholtz equation
\begin{eqnarray}\label{Helmholtz}
 \Delta u +k^2u=F \quad \mathrm{in} \ \Omega,
\end{eqnarray}
from boundary measurements $u|_\Gamma$ and $\partial_{\nu}u|_\Gamma$,
where $k>0$ is the wavenumber, $\Omega\subset \mathbb{R}^N$ ($N=2,3$) is a bounded Lipschitz domain with boundary $\Gamma$ and $\nu$ denotes the outward unit normal to $\Gamma$. A main difficulty of the inverse source problem with a single wavenumber is the non-uniqueness of the source due to the existence of non-radiating sources \cite{Albanese,Bleistein,Cheng2016, Dassios,Devaney1}, and several numerical methods with multi-frequency measurements \cite{Bao2015a, Bao2015b, Eller, Zhang} have been proposed to overcome it for the source with a compact support in the $L^2$ sense. However, fortunately, with a single wavenumber, the uniqueness can be obtained if {\it a priori} information
on the source is available \cite{El0,El2}. Readers may refer to \cite{El0,Badia,Badia1,Kang} and the reference therein for the further stability results.

In this paper, we assume that the source $F$ is a finite combination of well separated monopoles and dipoles of the form
\begin{equation}\label{FForm}
 F(x)=\sum_{j=1}^{M}(\lambda_j+\eta_j\cdot\nabla)\delta(x-z_j),
\end{equation}
where $\delta$ stands for the Dirac distribution, $M\in \mathbb{N}$ signifies the number of the source points, $\{z_j\}_{j=1}^M$ are points in $\Omega$,
and $\lambda_j$ and $\eta_j$ are respectively, scalar and vector source intensities such that
$|\lambda_j|+|\eta_j|\neq 0$ and $|\lambda_j\eta_j|=0, j=1,\cdots, M$. Here, the points $\{z_j\}_{j=1}^M$ are assumed to be mutually distinct.
The inverse source problem under concern is a location acquisition problem, which can be stated as: Find the locations $\{z_j\}_{j=1}^M$ of the source $F$ of the form \eref{FForm} in
the Helmholtz equation \eref{Helmholtz} from the boundary measurements
$u|_\Gamma$ and $\partial_{\nu}u|_\Gamma$ with a single wavenumber $k$.

Numerical methods for determining the multipolar sources have drawn a lot of interest in the literature.
For the Poisson equation, an algebraic method for recovering monopolar sources ($|\eta_j|=0$) was proposed in \cite{El1}, and has been extended to the case of multipolar sources in \cite{Chung1,Chung2,Nara1,Nara2}.
The algebraic method has also been developed to the 3D Helmholtz equation in \cite{El2} for monopolar sources,
and then to the 3D elliptic equations in \cite{Abdelaziz2} for multipolar sources. In the case of 2D elliptic equations, the method has been extended in \cite{Abdelaziz1} for monopolar sources. We refer to \cite{Kress} for a relevant paper on the reconstruction of extended sources for the 2D Helmholtz equation.

The purpose of this paper is to provide simple and effective numerical methods for determining
the multipolar sources for the Helmholtz equation in two and three dimensions.
We focus our attention on the non-iterative sampling-type methods, and we shall develop some novel direct sampling schemes in this paper.
The main technicality lies in the decaying property of oscillatory integrals, which is employed to construct indicator functions at any sampling point $z\in \Omega$
such that the proposed indicator functions attain the local maximum at $z_j$ in an open neighborhood of $z_j$.
We would like to emphasize that the proposed methods follow the one-shot ideas \cite{LLSS13,LLW13,LLZ14,L15} and possess the following attractive features: (a)
the indicator functions can be formulated in close form via the boundary data; (b) the imaging schemes are explicit since the imaging indicator functions do not rely on any matrix inversion or forward solution process;
(c) our methods are easy to implement with computational efficiency due to the fact that only cheap integrations are involved in the indicator functions;
(d) the locating schemes perform robust with respect to large noise level.

The outline of this paper is as follows. In Section~\ref{sec2}, we will present the rationale behind the new method
for the 2D and 3D case, respectively. In each case, rigorous mathematical justifications of the proposed sampling schemes are  provided in detail.
The direct sampling methods (DSM) and its two-level version are proposed at the end.
In Section~\ref{sec3}, uniqueness and stability  of the proposed methods are established and analyzed in terms of measurement error, respectively.
Numerical examples are demonstrated in Section~\ref{sec4} to verify the performance of the proposed reconstruction schemes.
Finally, we conclude the work with some remarks and discussions of future topics.

%%%%%%%%%%%%%%%%%%%%%%%%%%%%%%%%%%%%%%%%%%%%%%%%%%%
\section{Direct sampling method}\label{sec2}
%%%%%%%%%%%%%%%%%%%%%%%%%%%%%%%%%%%%%%%%%%%%%%%%%%%

In this section, we will present our sampling schemes and give rigorous mathematical justifications.
First, we make the following assumptions
\begin{eqnarray}\label{assumption1}
    L:=\min_{1\leq j, j'\leq M\atop j\neq j'} \mathrm{dist}(z_j,z_{j'})\gg \frac{2\pi}{k},
  \\[2mm]\label{assumption2}
   \frac{|\lambda_j|}{|\eta_{j'}|}\sim k  \quad  \mathrm{for}\ \lambda_j\neq0\ \mathrm{and}\  \eta_{j'}\neq0,
    \quad\forall j,j'=1,\cdots,M. ~~REMOVED!!
\end{eqnarray}
Condition \eref{assumption1} means that the sources are sparsely distributed, i.e., they are well separated
with respect to the wavelength $2\pi/k$. Condition \eref{assumption2} indicates that for larger wavenumber $k$, the magnitudes of the data produced by monopoles and dipoles should be of the same order. Otherwise, for example, if $|\lambda_j| \approx |\eta_{j'}|$, then the information from the monopoles would be annihilated by the data due to the dipoles,  such as
\[
\left|\frac{\e^{\i k |x|}}{4\pi |x|}\right|=\frac{1}{4\pi |x|}\quad \mathrm{and}\quad
\left|\nabla\frac{\e^{\i k |x|}}{4\pi |x|}\right|=\frac{k}{4\pi |x|}+\mathcal{O}(|x|^{-2});
\]
particularly, with the noisy data, one can not expect a reasonable numerical result on the monopoles.

%\begin{remark}
%Assumption \eref{assumption1} is needed only for the purpose of
%theoretical analysis. Our numerical examples in Section 4 show that the different
%sources could be relatively closer to each other. In particular, as
%long as the distance between different sources is bigger than a wavelength,
%the proposed imaging scheme would yield satisfactory reconstructions.
%\end{remark}

Let $\mathbb{S}^{N-1}:=\{x\in \mathbb{R}^N: |x|=1\}$ be the unit sphere in $\mathbb{R}^N$ and $\bold{d}\in\mathbb{S}^{N-1}$ be represented by
\[
d=\left\{ \begin{array}{ll}
(d_1, d_2), & N=2,\\
(d_1,d_2,d_3), & N=3,
\end{array}\right.
\]
For uniformity of exposition, we introduce a constant $d_0\equiv1$ throughout this paper. For any sampling point $z\in \Omega\subset\mathbb{R}^N$, we define $N+1$ indicator functions
\begin{equation}\label{indicator}
	I_{N,\ell}(z):=\frac{a_{N,\ell}}{ 2^{N-1}\pi}\int_{\mathbb{S}^{N-1}} \mathcal{R}(d) d_\ell \e^{-\i kd\cdot z}\d s(d),
	\quad \ell=0,1,\cdots, N,
\end{equation}
where
\begin{eqnarray}\label{coefficients}
	a_{N,\ell}:=\left\{ \begin{array}{ll}
		1, & \ell=0,\\
		\displaystyle\frac{N\i}{k}, & \ell=1,\cdots, N,
	\end{array}\right.
\end{eqnarray}
and
\begin{eqnarray}\label{RG}
\mathcal{R}(d):=\int_\Gamma \left(\e^{\i k x\cdot d}\partial_\nu u(x)- u(x)\partial_{\nu}\e^{\i k x\cdot d}\right) \d s(x).
\end{eqnarray}

These indicator functions will be used to find the locations $\{z_j\}_{j=1}^M$ of the sources.

\subsection{Two dimensional case}

To see the characteristics of the indicator functions in two dimensions, we need to establish two crucial lemmas.
\begin{lemma}\label{Asymptotic}%---------------------------------------------------------------------------
Let $z\in \mathbb{R}^2\backslash\{(0,0)\}$, $d=(d_1,d_2)\in \mathbb{S}^1$. Then
\begin{equation} \label{Asymptotic1}
\left|\int_{\mathbb{S}^1} d_p d_q \e^{\i k d\cdot z}\d s(d)\right|
= \mathcal{O} ( (k|z|)^{-1/2}), \quad k|z| \rightarrow \infty,
\end{equation}
where $p,q=0,1,2$.
\end{lemma}%--------------------------------------------------------------------------------------------

\Proof
Denote $z=|z|(\cos \alpha, \sin \alpha)$ and $d = (\cos \theta, \sin \theta)$,
then we have
\begin{eqnarray} \label{Integral1}
  \int_{\mathbb{S}^1}  \e^{\i k d\cdot z}\d s(d)
  =\int_0^{2\pi}  \e^{\i k|z|\cos(\theta-\alpha)}\d\theta.
\end{eqnarray}
Introduce the {\it Jacobi-Anger expansion} (see \cite[Section 2.22]{Watson} and \cite[p.75]{Colton})
%Introduce the {\it Jacobi-Anger expansion} (see \cite[Section 2.22]{Watson} and \cite[p.67]{Colton})
\begin{eqnarray} \label{Jacobi}
  \e^{\i \rho \cos \phi}=\sum_{n=-\infty}^{+\infty}\i^nJ_n(\rho)\e^{\i n \phi},
\end{eqnarray}
where $J_n$ is the Bessel function of the first kind of order $n$. By taking $\rho=k|z|$ and
$\phi=\theta-\alpha$ in \eref{Jacobi}, then integrating over $(0,2\pi)$ with respect to $\theta$
and using \eref{Integral1}, we obtain
\[
  \int_{\mathbb{S}^1}  \e^{\i k d\cdot z}\d s(d)=2\pi J_0(k|z|).
\]
Similarly, from $\cos \theta=(\e^{\i \theta}+\e^{-\i \theta})/2$,
$\sin\theta=(\e^{\i \theta}-\e^{-\i \theta})/{2\i}$, $\cos^2 \theta=(1+\cos2\theta)/2$,
$\sin^2\theta=(1-\cos2\theta)/2$ and \eref{Jacobi},  we derive
\begin{eqnarray}\nonumber
   \int_{\mathbb{S}^1} d_1 \e^{\i k d\cdot z}\d s(d)
  &=&\int_0^{2\pi} \e^{\i k|z|\cos(\theta-\alpha)} \cos \theta  \d\theta
 \\\nonumber
  &=&\sum_{n=-\infty}^{+\infty}\frac{\i^n}{2}J_n(k|z|)\int_0^{2\pi} \e^{\i n (\theta-\alpha)}
     (\e^{\i \theta}+\e^{-\i \theta})  \d\theta
 \\\nonumber
  &=&  \frac{\i^{-1}}{2}J_{-1}(k|z|)\e^{\i  \alpha}2\pi+ \frac{\i}{2}J_{1}(k|z|)\e^{-\i \alpha}2\pi
 \\\label{expansion1}
  &=& 2\pi \i\cos \alpha \ J_{1}(k|z|),
\end{eqnarray}
\[
  \int_{\mathbb{S}^1} d_2 \e^{\i k d\cdot z}\d s(d)
  = 2\pi \i\sin \alpha \ J_{1}(k|z|),
\]
\begin{eqnarray}\label{expansion2}
  \int_{\mathbb{S}^1} d_1^2 \e^{\i k d\cdot z}\d s(d)
  =\pi J_0(k|z|)- \pi \cos 2\alpha\ J_{2}(k|z|),
\end{eqnarray}
\[
   \int_{\mathbb{S}^1} d_2^2 \e^{\i k d\cdot z}\d s(d)
  =\pi J_0(k|z|)+ \pi \cos 2\alpha\ J_{2}(k|z|),
\]
and
\begin{eqnarray}\label{expansion3}
   \int_{\mathbb{S}^1} d_1d_2 \e^{\i k d\cdot z}\d s(d)
 =-\pi \sin 2\alpha\ J_{2}(k|z|).
\end{eqnarray}

Hence, the estimate \eref{Asymptotic1} follows from the asymptotic behavior (see \cite[p.74]{Colton})
%Further, by using the asymptotic behavior (see \cite[p.66]{Colton})
\[
   H^{(1)}_{n}(t)=\sqrt{\frac{2}{\pi t}}\e^{\i  (t-\frac{n\pi}{2}-\frac{\pi}{4})}
 \left\{1+\mathcal{O}(t^{-1})\right\}, \quad t\rightarrow \infty,
\]
and $J_n(t)= \mathrm{Re} \{ H^{(1)}_{n}(t)\}$, where $H^{(1)}_{n}$ denotes the Hankel function of the first kind of order $n$.  \qedhere
%-------------------------------------------------------------------------------------------------

The following lemma follows from the definition of the Bessel functions by truncating the
infinite series:
\begin{eqnarray*}
 J_0(t)&=&\sum_{p=0}^{\infty}\frac{(-1)^p}{(p!)^2}\left(\frac{t}{2}\right)^{2p},
\\
J_n(t)&=&\sum_{p=0}^{\infty}\frac{(-1)^p}{p!(n+p)!}\left(\frac{t}{2}\right)^{n+2p}
\\
&=&\frac{t^n}{2^nn!}\left(1+\sum_{p=1}^{\infty}(-1)^p\frac{\left(\frac{t}{2}\right)^{2p}}{p!(n+1)\cdots(n+p)}\right),
 \quad n\geq1.
\end{eqnarray*}
\begin{lemma}\label{besselj}%---------------------------------------------------------------------------
For $0<t<1$, we have
\begin{eqnarray}
   0<J_{0}(t)< 1-\frac{t^2}{4}+\frac{t^4}{64},\label{besselj0} \\
   0<J_{1}(t)<\frac{t}{2},\label{besselj1} \\
   0<J_{2}(t)<\frac{t^2}{8}.\label{besselj2}
\end{eqnarray}
\end{lemma}%--------------------------------------------------------------------------------------------

%---------------------------------------------------------------------------------------------

Now, we present the indicating behaviors of $\{I_{2,\ell}(z)\}_{\ell=0}^2$,
which play a vital role in our schemes.

\begin{theorem}\label{behavior_2D}
	Let source $F$ be of the form \eref{FForm} with $\eta_j=(\eta_{j,1},\eta_{j,2})$ satisfying $|\lambda_j|+|\eta_j|\neq 0$ and $|\lambda_j\eta_j|=0$, the assumptions
	\eref{assumption1}, \eref{assumption2} hold, and $\{I_{2,\ell}(z)\}_{\ell=0}^2$ be described as in \eref{indicator}. Then, we have the following asymptotic expansions
	\begin{eqnarray}
	& I_{2,0}(z_j) = \lambda_j+\mathcal{O}\left((kL)^{-1/2}\right), \label{behavior_2D_1} \\
	& I_{2,\ell}(z_j) = \eta_{j,\ell}+\mathcal{O}\left((kL)^{-1/2}\right), \quad \ell=1,2. \label{behavior_2D_2}
	\end{eqnarray}
	Moreover, there exists an open neighborhood of
	$z_j$, $neigh(z_j)$, such that
	\begin{eqnarray*}
	& |I_{2,0}(z)| \leq |\lambda_j| +\mathcal{O}\left((kL)^{-1/2}\right),\quad z\in  neigh(z_j) \quad for \; \lambda_j  \neq 0,\\
	& |I_{2,\ell}(z)| \leq |\eta_{j,\ell}|+\mathcal{O}\left((kL)^{-1/2}\right),\quad z\in  neigh(z_j) \quad for \; \eta_{j,\ell} \neq 0,\quad \ell=1,2,
	\end{eqnarray*}
	where the equalities hold only at $z=z_j$. That is, $z_j$ is a local maximizer of $|I_{2,0}(z)|$ for $\lambda_j  \neq 0$
	and $|I_{2,\ell}(z)|$ for $\eta_{j,\ell} \neq 0$ in $neigh(z_j)$.
\end{theorem}

\Proof
Without loss of generality, we only consider the indicating behavior of $I_{2,1}(z)$.
First, by multiplying equation \eref{Helmholtz} by  $\e^{\i k x\cdot d}$ for $d\in \mathbb{S}^1$, then integrating over $\Omega$,
using Green's formula and \eref{RG} we obtain
\begin{eqnarray}\label{sum_2D}
  \sum_{j=1}^{M}\lambda_j\e^{\i k d\cdot z_j}- \i k \sum_{j=1}^{M}(\eta_j\cdot d) \e^{\i k d\cdot z_j}=\mathcal{R}(d),\quad d \in  \mathbb{S}^1.
\end{eqnarray}
Next, we multiply equation \eref{sum_2D} by $d_1 \e^{-\i k d\cdot z}$ for sampling point $z\in \Omega$, integrate over $\mathbb{S}^1$, and then derive
\begin{eqnarray*}
    I_{2,1}(z)&=&\frac{1}{\pi}\sum_{j=1}^{M}\int_{\mathbb{S}^1}\left(\eta_{j,1}d_1d_1
      +\eta_{j,2}d_1 d_2\right) \e^{\i k d\cdot (z_j-z)}\d s(d) \\
   &&+\frac{\i}{k\pi}\sum_{j=1}^{M}\lambda_j\int_{\mathbb{S}^1}d_1 \e^{\i k d\cdot (z_j-z)} \d s(d).
\end{eqnarray*}
Further, by using the assumptions \eref{assumption1}, \eref{assumption2}, and Lemma \ref{Asymptotic}, we have for $j'=1,...,M$,
\begin{eqnarray*}
    I_{2,1}(z_{j'})&=&\eta_{j',1}
   +\frac{1}{\pi}\sum_{j=1\atop j\neq j'}^{M}\int_{\mathbb{S}^1}\left(\eta_{j,1}d_1d_1
      +\eta_{j,2}d_1 d_2\right) \e^{\i k d\cdot (z_j-z_{j'})}\d s(d)
 \\
  && +\frac{\i}{k\pi}\sum_{j=1\atop j\neq j'}^{M}\lambda_j\int_{\mathbb{S}^1}d_1 \e^{\i k d\cdot (z_j-z_{j'})} \d s(d)
 \\
  &=&  \eta_{j',1} +\mathcal{O}\left((kL)^{-1/2}\right),
\end{eqnarray*}
which leads to the expansion \eref{behavior_2D_2}.

Let $B_\rho(z_{j'}):=\{z: |z-z_{j'}|<\rho\}$ and $0<\rho < k^{-1}(kL)^{-1/2}$. Then we have
\[
    |z_j-z|\geq L-\rho\gg \frac{2\pi}{k},\quad z\in B_\rho(z_{j'}), \ j\neq j',
\]
which, together with Lemma \ref{Asymptotic}, yields
\begin{eqnarray*}
    I_{2,1}(z)&=& \frac{1}{\pi}\int_{\mathbb{S}^1}\left(\eta_{j',1}d_1d_1
      +\eta_{j',2}d_1 d_2\right) \e^{\i k d\cdot (z_{j'}-z)}\d s(d) \\
   &&+\frac{\i}{k\pi}\lambda_{j'}\int_{\mathbb{S}^1}d_1 \e^{\i k d\cdot (z_{j'}-z)} \d s(d)
    +\mathcal{O}\left((kL)^{-1/2}\right), \quad z\in B_\rho(z_{j'}).
\end{eqnarray*}
By using \eref{expansion1}--\eref{besselj2}, we deduce that for $z\in B_\rho(z_{j'})$ and $z\neq z_{j'}$,
\begin{eqnarray*}
     \frac{1}{\pi}\left|\int_{\mathbb{S}^1}d_1d_1 \e^{\i k d\cdot (z_{j'}-z)}\d s(d)\right|
   &\leq &J_0(k|z_{j'}-z|)+ J_{2}(k|z_{j'}-z|) \\
    &< & 1-\frac{\tau^2}{8}+\frac{\tau^4}{64} \\
    &< & 1,
\end{eqnarray*}
\[
   \frac{1}{\pi}\left|\int_{\mathbb{S}^1}d_1d_2 \e^{\i k d\cdot (z_{j'}-z)}\d s(d)\right|
   \leq J_{2}(k|z_{j'}-z|)<\frac{\tau^2}{8}
\]
and
\[
   \frac{1}{k\pi}\left|\int_{\mathbb{S}^1}d_1 \e^{\i k d\cdot (z_{j'}-z)}\d s(d)\right|
   \leq \frac{2}{k}J_{1}(k|z_{j'}-z|)<\frac{\tau}{k},
\]
where $\tau= k|z_{j'}-z|$. Hence, we obtain
\[
  |I_{2,1}(z)|< |\eta_{j',1}|+\mathcal{O}((kL)^{-1/2})
  \quad \mathrm{for}\  z\in B_\rho(z_{j'})\backslash\{z_{j'}\}.
\]
This completes the proof of the theorem.  \qedhere
%------------------------------------------------------------------------------------------------

%%%%%%%%%%%%%%%%%%%%%%%%%%%%%%%%%%%%%%%%%%%%%%%%%
\subsection{Three dimension case}
%%%%%%%%%%%%%%%%%%%%%%%%%%%%%%%%%%%%%%%%%%%%%%%%%

In this subsection, we shall analyze the properties of the indicator functions in 3D. Before showing the indicating behaviors, we establish two crucial lemmas.

\begin{lemma}\label{lem:3dAsymptotic}%-----------------------------------------------------------------------

Let $z\in \mathbb{R}^3\backslash\{(0,0,0)\}$, $d=(d_1,d_2,d_3)\in \mathbb{S}^2$. Then
\begin{eqnarray} \label{3dAsymptotic}
\left|\int_{\mathbb{S}^2} d_p d_q \e^{\i k d\cdot z}\d s(d)\right|
= \mathcal{O} ( (k|z|)^{-1}), \quad k|z| \rightarrow \infty,
\end{eqnarray}
where $p,q=0, 1,2,3$.
\end{lemma}%--------------------------------------------------------------------------------------------

The proof is technical and lengthy, therefore we deter it till Appendix A.

%--------------------------------------------------------------------------------------------------------
The following lemma is a 3D version of Lemma~\ref{besselj} and follows from the definition of
the  spherical Bessel functions
\[
j_n(t)=\sum_{p=0}^{\infty}\frac{(-1)^pt^{n+2p}}{2^p p!1\cdot3\cdot\cdot\cdot(2n+2p+1)}.
\]
\begin{lemma}\label{3dbesselj}%---------------------------------------------------------------------------
For $0<t<1$, we have
\begin{eqnarray} \label{3dbesselj0}
  0<j_{0}(t)< 1-\frac{t^2}{6}+\frac{t^4}{120},
  \\\label{3dbesselj1}
   0<j_{1}(t)<\frac{t}{3},
\\\label{3dbesselj2}
   0<j_{2}(t)<\frac{t^2}{15}.
\end{eqnarray}
\end{lemma}%--------------------------------------------------------------------------------------------

%-------------------------------------------------------------------------------------------------------------

Now, the indicating behaviors of $\{I_{3,\ell}(z)\}_{\ell=0}^3$ are presented in the following theorem.

\begin{theorem}\label{behavior_3D}
Let source $F$ be of the form \eref{FForm} with $\eta_j=(\eta_{j,1},\eta_{j,2},\eta_{j,3})$ satisfying $|\lambda_j|+|\eta_j|\neq 0$ and $|\lambda_j\eta_j|=0$, the assumptions
\eref{assumption1}, \eref{assumption2} hold, and indicator functions $\{I_{3,\ell}(z)\}_{\ell=0}^3$ be described as in \eref{indicator}. Then, we have the following asymptotic expansions
\begin{eqnarray}
& I_{3,0}(z_j) = \lambda_j+\mathcal{O}\left((kL)^{-1}\right), \label{behavior_3D_1} \\
& I_{3,\ell}(z_j) = \eta_{j,\ell}+\mathcal{O}\left((kL)^{-1}\right),\quad \ell= 1,2,3. \label{behavior_3D_2}
\end{eqnarray}
Moreover, there exists an open neighborhood of $z_j$, $neigh(z_j)$, such that
\begin{eqnarray*}
& |I_{3,0}(z)|\leq |\lambda_j| +\mathcal{O}\left((kL)^{-1}\right),\quad z\in  neigh(z_j) \quad for \; \lambda_j  \neq 0, \\
& |I_{3,\ell}(z)| \leq |\eta_{j,\ell}|+\mathcal{O}\left((kL)^{-1}\right),\quad z\in  neigh(z_j) \quad for \; \eta_{j,\ell} \neq 0,\quad\ell=1,2,3,
\end{eqnarray*}
where the equalities hold only at $z=z_j$. That is, $z_j$ is a local maximizer of $|I_{3,0}(z)|$ for $\lambda_j \neq 0$
and $|I_{3,\ell}(z)|$ for $\eta_{j,\ell} \neq 0$ in $neigh(z_j)$.
\end{theorem}
\Proof
Without loss of generality, we only consider the indicating behavior of $I_{3,1}(z)$.
First, by multiplying equation \eref{Helmholtz} by  $\e^{\i k x\cdot d}$ for $d\in \mathbb{S}^2$, then integrating over $\Omega$, using Green's formula and \eref{RG} we obtain
\begin{eqnarray}\label{sum}
  \sum_{j=1}^{M}\lambda_j\e^{\i k d\cdot z_j}- \i k \sum_{j=1}^{M}(\eta_j\cdot d) \e^{\i k d\cdot z_j}=\mathcal{R}(d),
  \quad d \in  \mathbb{S}^2.
\end{eqnarray}
Next, we multiply equation \eref{sum} by $d_1 \e^{-\i k d\cdot z}$ for sampling point $z\in \Omega$, integrate over $\mathbb{S}^2$,
and then derive
\begin{eqnarray*}
    I_{3,1}(z)&=&\frac{3}{4\pi}\sum_{j=1}^{M}\int_{\mathbb{S}^2}\left(\eta_{j,1}d_1d_1
      +\eta_{j,2}d_1 d_2+\eta_{j,3}d_1 d_3\right) \e^{\i k d\cdot (z_j-z)}\d s(d)
 \\
   &&+\frac{3\i}{4k\pi}\sum_{j=1}^{M}\lambda_j\int_{\mathbb{S}^2}d_1 \e^{\i k d\cdot (z_j-z)} \d s(d).
\end{eqnarray*}
Further, by using the assumptions \eref{assumption1}, \eref{assumption2}, and Lemma \ref{lem:3dAsymptotic}, we have for $j'=1,\cdots,M$,
\begin{eqnarray*}
    I_{3,1}(z_{j'})&=&\eta_{j',1}
   +\frac{3}{4\pi}\sum_{j=1\atop j\neq j'}^{M}\int_{\mathbb{S}^2}\left(\eta_{j,1}d_1d_1
      +\eta_{j,2}d_1 d_2+\eta_{j,3}d_1 d_3\right) \e^{\i k d\cdot (z_j-z_{j'})}\d s(d)
 \\
  && +\frac{3\i}{4k\pi}\sum_{j=1\atop j\neq j'}^{M}\lambda_j\int_{\mathbb{S}^2}d_1 \e^{\i k d\cdot (z_j-z_{j'})} \d s(d)
 \\
  &=&  \eta_{j',1} +\mathcal{O}\left((kL)^{-1}\right),
\end{eqnarray*}
which leads to the expansion \eref{behavior_3D_2}.

Let $B_\rho(z_{j'}):=\{z: |z-z_{j'}|<\rho\}$ and $0<\rho < k^{-1}(kL)^{-1}$. Then we have
\begin{eqnarray*}
    |z_j-z|\geq L-\rho\gg \frac{2\pi}{k},\quad z\in B_\rho(z_{j'}), \ j\neq j',
\end{eqnarray*}
which, together with Lemma \ref{lem:3dAsymptotic}, yields
\begin{eqnarray*}
    I_{3,1}(z)&=& \frac{3}{4\pi}\int_{\mathbb{S}^2}\left(\eta_{j',1}d_1d_1
      +\eta_{j',2}d_1 d_2+\eta_{j',3}d_1 d_3\right) \e^{\i k d\cdot (z_{j'}-z)}\d s(d)
  \\
   &&+\frac{3\i}{4k\pi}\lambda_{j'}\int_{\mathbb{S}^2}d_1 \e^{\i k d\cdot (z_{j'}-z)} \d s(d)
    +\mathcal{O}\left((kL)^{-1}\right),
   \quad z\in B_\rho(z_{j'}).
\end{eqnarray*}
By using \eref{3d0}--\eref{3dbesselj2}, we deduce that for $z\in B_\rho(z_{j'})$ and $z\neq z_{j'}$,
\begin{eqnarray*}
     \frac{3}{4\pi}\left|\int_{\mathbb{S}^2}d_1d_1 \e^{\i k d\cdot (z_{j'}-z)}\d s(d)\right|
   &\leq &j_0(k|z_{j'}-z|)
  \\
  &&+ \frac{|3\cos^2\alpha-1|+\sin^2\alpha}{2}j_{2}(k|z_{j'}-z|)
   \\
  &\leq &j_0(k|z_{j'}-z|)+j_{2}(k|z_{j'}-z|)
 \\
    &< & 1-\frac{\tau^2}{10}+\frac{\tau^4}{120}
 \\
    &< & 1,
\end{eqnarray*}
\begin{eqnarray*}
     \frac{3}{4\pi}\left|\int_{\mathbb{S}^2}d_1d_2 \e^{\i k d\cdot (z_{j'}-z)}\d s(d)\right|
   \leq  \frac{3}{2}j_{2}(k|z_{j'}-z|)<\frac{\tau^2}{10},
\end{eqnarray*}
\begin{eqnarray*}
     \frac{3}{4\pi}\left|\int_{\mathbb{S}^2}d_1d_3 \e^{\i k d\cdot (z_{j'}-z)}\d s(d)\right|
   \leq  3j_{2}(k|z_{j'}-z|)<\frac{\tau^2}{5}
\end{eqnarray*}
and
\begin{eqnarray*}
     \frac{3}{4k\pi}\left|\int_{\mathbb{S}^2}d_1 \e^{\i k d\cdot (z_{j'}-z)}\d s(d)\right|
   \leq \frac{3}{k}j_{1}(k|z_{j'}-z|)<\frac{\tau}{k},
\end{eqnarray*}
where $\tau= k|z_{j'}-z|$. Hence, we obtain
\begin{eqnarray*}
    |I_{3,1}(z)|< |\eta_{j',1}|+\mathcal{O}\left((kL)^{-1}\right)
  \quad \mathrm{for}\  z\in B_\rho(z_{j'})\ \mathrm{and}\ z\neq z_{j'}.
\end{eqnarray*}
This completes the proof of the theorem.  \qedhere
%----------------------------------------------------------------------------------------------------

\subsection{Sampling schemes}

Now we are in the position to present the main theorem of this paper by combining Theorems~\ref{behavior_2D} and \ref{behavior_3D}.
\begin{theorem}\label{indicating_behavior}
	Let source $F$ be of the form \eref{FForm} with $\eta_j=(\eta_{j,1},\cdots,\eta_{j,N})$ satisfying $|\lambda_j|+|\eta_j|\neq 0$ and $|\lambda_j\eta_j|=0$, the assumptions
	\eref{assumption1}, \eref{assumption2} hold, and indicator functions $\{I_{N,\ell}(z)\}_{\ell=0}^N$ be described as in \eref{indicator}. Then, we have the following asymptotic expansions
	\begin{eqnarray}
	& I_{N,0}(z_j) = \lambda_j+\mathcal{O}\left((kL)^{-\frac{N-1}{2}}\right),\label{behavior_1}\\
	& I_{N,\ell}(z_j) = \eta_{j,\ell}+\mathcal{O}\left((kL)^{-\frac{N-1}{2}}\right), \quad \ell=1,\cdots, N. \label{behavior_2}
	\end{eqnarray}
	Moreover, there exists an open neighborhood of $z_j$, $neigh(z_j)$, such that for all $z\in  neigh(z_j)$, it holds that
	\begin{eqnarray*}
	|I_{N,0}(z)|\leq |\lambda_j| +\mathcal{O}\left((kL)^{-\frac{N-1}{2}}\right)\quad for \; \lambda_j  \neq 0, \\
	|I_{N,\ell}(z)| \leq |\eta_{j,\ell}|+\mathcal{O}\left((kL)^{-\frac{N-1}{2}}\right)\quad for \; \eta_{j,\ell} \neq 0,\quad \ell=1,\cdots, N,
	\end{eqnarray*}
	where the equalities hold only at $z=z_j$. That is, $z_j$ is a local maximizer of $|I_{N,0}(z)|$ for $\lambda_j \neq 0$ and $|I_{N,\ell}(z)|$ for $\eta_{j,\ell} \neq 0$ in $neigh(z_j)$.
\end{theorem}

Motivated by Theorem \ref{indicating_behavior}, we are ready to present the first direct sampling method in $\mathbb{R}^N$, see Algorithm DSM.
%\begin{itemize}
%	\item Collect the Cauchy data $u|_\Gamma$ and $\partial_\nu u|_\Gamma$ on the boundary $\Gamma$ with the fixed wavenumber $k$;
%	\item Select a global coarse sampling grid $\tau$ over the probe region $\Omega$;
%	\item For each sampling point $z\in \tau$, compute the values of the indicator functions
%	$\{I_{N,l}(z)\}_{l=0}^N$;
%	\item According to the image of $|I_{N,l}(z)|, z\in\tau, l=0,1,\cdots, N$, collect all the significant local maximizers $\tilde{z}_{j,l}, j=1,\cdots,\tilde{M},l=0,1,\cdots,N$, where $\tilde{M}$ denotes the estimated number of sources. If several significant local maximizers are clustered in a region whose diameter is far less than $4\pi/k$, then they are treated as a single local maximizer;
%	\item For each local maximizer $\tilde{z}_{j,l}, j=1,\cdots,\tilde{M},l=0,1,\cdots,N$, select a local fine sampling grid $\tau_{j,l}$ centered at $\tilde{z}_{j,l}$ such that its diameter is less than $4\pi/k$;
%	\item Locate the maximizers $z_{j,l}$ of $|I_{N,l}(z)|$ over $\tau_{j,l}$ for $j=1,\cdots,\tilde{M}, l=0,1,\cdots, N$;
%	\item Categorize $z_{j,l}, j=1,\cdots,\tilde{M}, l=0,1,\cdots,N$, into $\tilde{M}$ clustered groups such that in each group the distances between distinct points are less than $2\pi/k$;
%	\item Take the average location of all the points in each group as the reconstruction of the
%	locations of the sources;
%\end{itemize}

\begin{table}[h]
\centering
\begin{tabular}{cp{.8\textwidth}}
\toprule
\multicolumn{2}{l}{{\bf Algorithm DSM:}\quad Direct sampling method for recovering point sources} \\
\midrule
 {\bf Step 1} & Collect the Cauchy data $u|_\Gamma$ and $\partial_\nu u|_\Gamma$  with the fixed wavenumber $k$;  \\
{\bf Step 2} & Select a sampling grid $\tau$ over the probe region $\Omega$; \\
{\bf Step 3} & For each sampling point $z\in \tau$, compute the values of the indicator functions
	$\{I_{N,\ell}(z)\}_{\ell=0}^N$; \\
{\bf Step 4} & According to the images of $|I_{N,\ell}(z)|, z\in\tau, \ell=0,1,\cdots, N$, collect all the significant local maximizers $z_{j,\ell}, j=1,\cdots,\tilde{M}, \ell=0,1,\cdots,N$, where $\tilde{M}$ denotes the estimated number of sources. If several significant local maximizers are clustered in a region whose diameter is far less than $4\pi/k$, then they are treated as a single local maximizer; \\
{\bf Step 5} & Cluster $z_{j,\ell}, j=1,\cdots,\tilde{M}, \ell=0,1,\cdots,N$, into $\tilde{M}$  groups such that in each group the distances between distinct points are less than $2\pi/k$; \\
{\bf Step 6} & Take the average location of all the points in each group as the reconstruction of the
	locations of the sources; \\
\bottomrule
\end{tabular}
\end{table}

\begin{figure}
	\centering
	\subfloat[]{\includegraphics[width=0.45\textwidth]{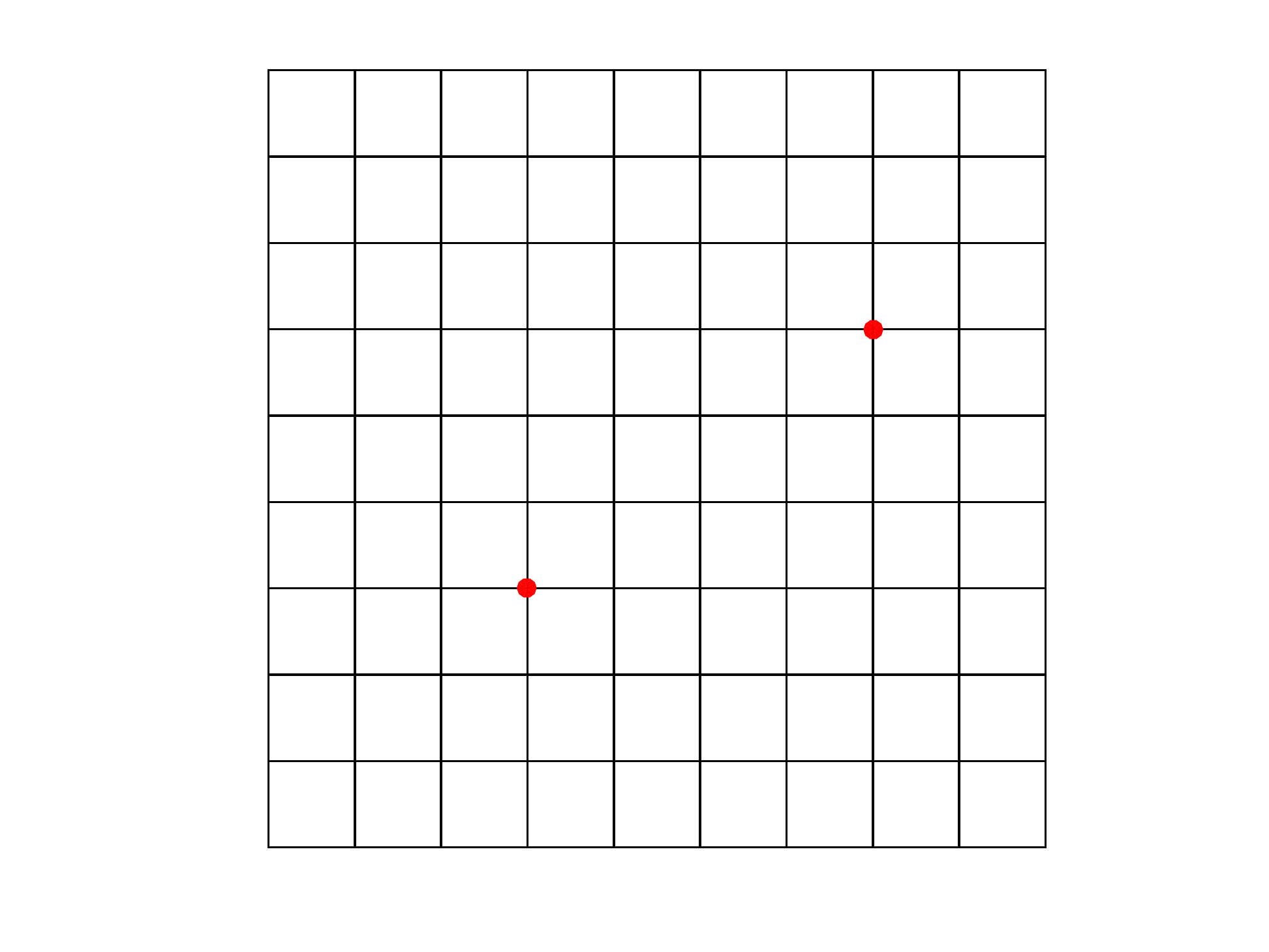}}
	\subfloat[]{\includegraphics[width=0.45\textwidth]{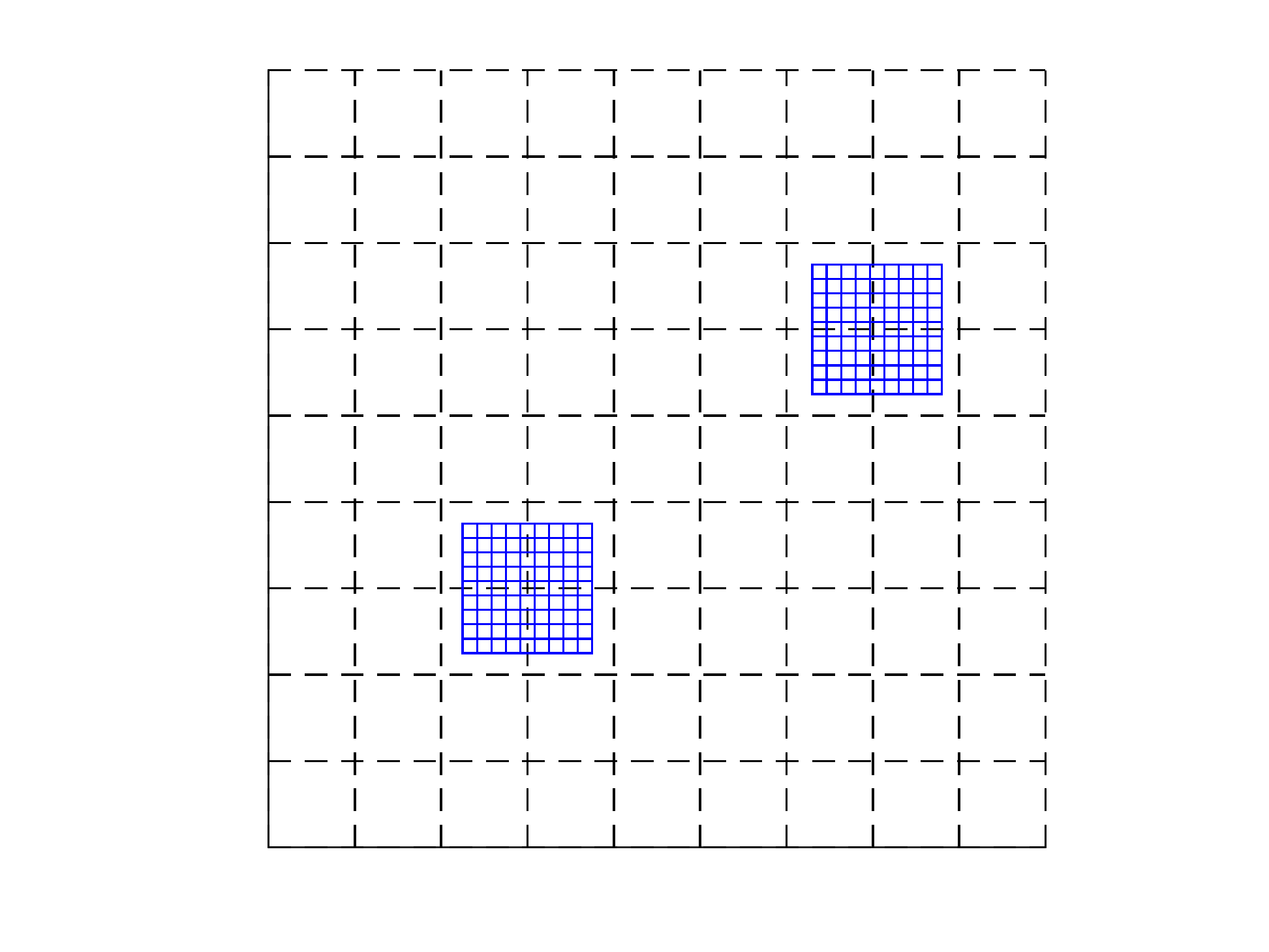}}
	\caption{An illustration of DSM2 in 2D. (a) global coarse grid (black lines). The maximizers of the indicator function over the global coarse grid are denoted by the small red points; (b) local fine grid (blue lines). In the vicinity of each maximizer, a local fine grid is used to resample the concerned region and enhance the accuracy of the maximizer location.}
	\label{fig:DSM2}
\end{figure}

Our numerical experiments showed that the sampling scheme DSM could work very well if the sampling grid is sufficiently fine.
However, a uniform fine mesh is computationally expensive and it would be more efficient to use a relatively coarse sampling grid for the most part of the probe region $\Omega$ since the point sources are assumed to be sparsely distributed. On the other hand, it is difficult to locate the local maximizers accurately if only a uniformly coarse grid is available. To guarantee the accuracy and computational efficiency simultaneously, we propose a two-level sampling strategy, whose main feature consists of a global coarse grid for rough locating and a local fine grid for fine tuning. The two-level direct sampling scheme is stated as Algorithm DSM2.
To sketch the idea of DSM2,  we refer to Figure~\ref{fig:DSM2}.

\begin{table}[h]
\centering
\begin{tabular}{cp{.8\textwidth}}
\toprule
\multicolumn{2}{l}{{\bf Algorithm DSM2:}\quad Two-level DSM for recovering point sources} \\
\midrule
 {\bf Step 1} & Collect the Cauchy data $u|_\Gamma$ and $\partial_\nu u|_\Gamma$  with the fixed wavenumber $k$;  \\
{\bf Step 2} & Select a global coarse sampling grid $\tau$ over the probe region $\Omega$; \\
{\bf Step 3} & For each sampling point $z\in \tau$, compute the values of the indicator functions
	$\{I_{N,\ell}(z)\}_{\ell=0}^N$; \\
{\bf Step 4} & According to the images of $|I_{N,\ell}(z)|, z\in\tau, \ell=0,1,\cdots, N$, collect all the significant local maximizers $\tilde{z}_{j,\ell}, j=1,\cdots,\tilde{M}, \ell=0,1,\cdots,N$, where $\tilde{M}$ denotes the estimated number of sources. If several significant local maximizers are clustered in a region whose diameter is far less than $4\pi/k$, then they are treated as a single local maximizer; \\
{\bf Step 5} & For each local maximizer $\tilde{z}_{j,\ell}, j=1,\cdots,\tilde{M}, \ell=0,1,\cdots,N$, select a local fine sampling grid $\tau_{j,\ell}$ centered at $\tilde{z}_{j,\ell}$ such that its diameter is less than $4\pi/k$;  \\
{\bf Step 6} & Locate the maximizers $z_{j,\ell}$ of $|I_{N,\ell}(z)|$ over $\tau_{j,\ell}$ for $j=1,\cdots,\tilde{M}, \ell=0,1,\cdots, N$; \\
{\bf Step 7} & Cluster $z_{j,\ell}, j=1,\cdots,\tilde{M}, \ell=0,1,\cdots,N$, into $\tilde{M}$  groups such that in each group the distances between distinct points are less than $2\pi/k$; \\
{\bf Step 8} & Take the average location of all the points in each group as the reconstruction of the
	locations of the sources; \\
\bottomrule
\end{tabular}
\end{table}

\begin{remark}
In addition, if we know a priori that $\eta_j=0,\forall j$ (resp. $\lambda_j=0,\forall j$), i.e., the target source consists of only monopoles (resp. dipoles), then  according to Theorem \ref{indicating_behavior}, the above algorithms could be simplified by utilizing only indicator(s) $I_{N,0}$ (resp. $I_{N,\ell}, \ell=1,\cdots,N$).
\end{remark}

%%%%%%%%%%%%%%%%%%%%%%%%%%%%%%%%%%%%%%%%%%%%%%%%%%
\section{Uniqueness and stability}\label{sec3}
%%%%%%%%%%%%%%%%%%%%%%%%%%%%%%%%%%%%%%%%%%%%%%%%%%

We begin this section with an interesting result on the uniqueness under the assumptions
\eref{assumption1} and \eref{assumption2},
which indicates that our sampling schemes can well distinguish the sparsely distributed sources.

\begin{theorem}\label{Unique}
Let source $F$ be of the form \eref{FForm} and the assumptions \eref{assumption1}, \eref{assumption2} hold, then the source  $F$ can be uniquely determined by the boundary measurements $u|_\Gamma$ and $\partial_\nu u|_\Gamma$.
\end{theorem}

\Proof
Without loss of generality, we only need to consider the homogeneous boundary value problem for the 3D case. Let
$u|_\Gamma=\partial_\nu u|_\Gamma=0$. Then we have $\mathcal{R}(d)=0$, $\forall d \in \mathbb{S}^{N-1}$, and thus, $I_{N,\ell}(z)=0,\forall z\in \Omega, \ell=0,1,\cdots,N$. In particular, $I_{N,\ell}(z_j)=0, \ell=0,1,\cdots,N$. Furthermore, by using \eref{behavior_1} and \eref{behavior_2}, we
obtain
\[
\lambda_j=\mathcal{O}\left((kL)^{-\frac{N-1}{2}}\right)\ \mathrm{and}\
\eta_{j,\ell}=\mathcal{O}\left((kL)^{-\frac{N-1}{2}}\right),\  \ell=1,\cdots, N,
\]
which, together with assumptions \eref{assumption1} and \eref{assumption2},
yields $\lambda_j=0$ and $\eta_{j,\ell}=0, \ell=1,\cdots, N$. This completes the proof.   \qedhere

In the following, we are going to analyze the stability of our methods.
Let the measured noisy data $u^\epsilon, \partial_\nu u^\epsilon\in L^2(\Gamma)$ satisfy
\begin{eqnarray}\label{noise}
    \|u^\epsilon-u\|_{L^2(\Gamma)} \leq \epsilon \|u\|_{L^2(\Gamma)}, \qquad
    \|\partial_\nu u^\epsilon-\partial_\nu u\|_{L^2(\Gamma)} \leq \epsilon\|\partial_\nu u\|_{L^2(\Gamma)},
\end{eqnarray}
where $0<\epsilon\ll1$.

We introduce the perturbed indicator functions
\begin{equation}\label{indicatorP}
I_{N,\ell}^\epsilon(z):=\frac{a_{N,\ell}}{2^{N-1}\pi}\int_{\mathbb{S}^{N-1}} \mathcal{R}^\epsilon(d) d_\ell \e^{-\i kd\cdot z}\d s(d), \quad \ell=0,1,\cdots, N,
\end{equation}
where $\{a_{N,\ell}\}_{\ell=0}^N$ are defined by \eref{coefficients} and
\begin{eqnarray}\label{perturbed_RG}
\mathcal{R}^\epsilon(d):=\int_\Gamma \left(\e^{\i k x\cdot d}\partial_\nu u^\epsilon(x)- u^\epsilon(x)\partial_{\nu}\e^{\i k x\cdot d}\right) \d s(x), \quad d \in \mathbb{S}^N.
\end{eqnarray}
Then, we have

\begin{theorem}\label{stability}
Let source $F$ be of the form \eref{FForm}  with $\eta_j=(\eta_{j,1},\cdots,\eta_{j,N})$ satisfying $|\lambda_j|+|\eta_j|\neq 0$ and $|\lambda_j\eta_j|=0$, the assumptions
\eref{assumption1}, \eref{assumption2} hold, and indicator functions $\{I_{N,\ell}^\epsilon(z)\}_{\ell=0}^N$ be described as in \eref{indicatorP}. Then, we have the following
asymptotic expansions
\begin{eqnarray}
   I^\epsilon_{N,0}(z_j) = \lambda_j+\mathcal{O}\left((kL)^{-\frac{N-1}{2}}\right)+\mathcal{O}(\epsilon),
   \label{Pbehavior0}\\
   I^\epsilon_{N,\ell}(z_j) = \eta_{j,\ell}+\mathcal{O}\left((kL)^{-\frac{N-1}{2}}\right)+\mathcal{O}(\epsilon), \quad\ell=1,\cdots,N. \label{Pbehaviorl}
\end{eqnarray}
Moreover, there exists an open neighborhood of $z_j$, $neigh(z_j)$, such that for all $z\in  neigh(z_j)$, it holds that
\begin{eqnarray}
   |I^\epsilon_{N, 0}(z)| \leq |\lambda_j| +\mathcal{O}\left((kL)^{-\frac{N-1}{2}}\right)+\mathcal{O}(\epsilon)\ for \; \lambda_j  \neq 0, \label{Pbehavior2}\\
   |I^\epsilon_{N,\ell}(z)| \leq |\eta_{j,\ell}|+\mathcal{O}\left((kL)^{-\frac{N-1}{2}}\right)+\mathcal{O}(\epsilon)\ for \; \eta_{j,\ell} \neq 0,\ \ell=1,\cdots,N, \label{Pbehavior3}
\end{eqnarray}
where the equalities hold only at $z=z_j$. That is, $z_j$ is a local maximizer of $|I^\epsilon_{N,0}(z)|$ for $\lambda_j \neq 0$ and $|I^\epsilon_{N,\ell}(z)|$ for $\eta_{j,\ell} \neq 0$ in $neigh(z_j)$.
\end{theorem}

\Proof
From \eref{noise} and the definition of $\mathcal{R}^\epsilon$ in \eref{perturbed_RG}, it can be readily seen that
\[
    \left\|\mathcal{R}^\epsilon-\mathcal{R}\right\|_{L^2(\mathbb{S}^N)} =\mathcal{O}(\epsilon).
\]
Then, by a similar argument as in the proofs of Theorem \ref{behavior_2D} and Theorem \ref{behavior_3D}, we complete the proof of the theorem. \qedhere

\begin{remark}\label{remark_stability}
Theorem \ref{stability} shows that small perturbation $\epsilon$ would not essentially affect the recovery of the locations. Actually, the indicator functions are also robust with respect to large noise. For example, the representation of
\[
I_{N,0}(z)=\frac{1}{2^{N-1} \pi}\int_{\mathbb{S}^{N-1}}   \e^{-\i kd\cdot z}
\int_\Gamma \left(\e^{\i k x\cdot d}\partial_\nu u(x)- u(x)\partial_{\nu}\e^{\i k x\cdot d}\right) \d s(x)\d s(d)
\]
possesses two generalized Fourier integrals, which could significantly filter the effect of the noise. This indicates that for our sampling schemes,  the reconstruction of locations $z_j$ would be inherently insensitive to the measurement noise. The robustness will be shown in our numerical examples in Section \ref{sec4}.
\end{remark}

%%%%%%%%%%%%%%%%%%%%%%%%%%%%%%%%%%%%%%%%%%%%%%%%%%
\section{Numerical examples}\label{sec4}
%%%%%%%%%%%%%%%%%%%%%%%%%%%%%%%%%%%%%%%%%%%%%%%%%%

In this section, we present extensive numerical results to demonstrate
the feasibility and effectiveness of the proposed methods. We will specify details of the numerical implementation of the DSM and DSM2. If unspecified elsewhere, the reconstructed coordinates of the    source points were obtained by the scheme DSM2 in the experiments.

\subsection{Synthetic data}

Synthetic Cauchy data was generated by using the closed form of the radiating fields. Recall the fundamental solution to the Helmholtz equation
\[
\Phi(x; z)=
\left\{\begin{array}{ll}
\displaystyle\frac{\i}{4}H_0^{(1)}(k|x-z|), & x\in\mathbb{R}^2, \\
\vspace{-3mm}&\\
\displaystyle\frac{\e^{\i k |x-z|}}{4\pi |x-z|}, & x\in\mathbb{R}^3, \\
\end{array}\right.
\]
where $H^{(1)}_{0}$ is the Hankel function of the
first kind of order zero. Then the unique solution to the Helmholtz equation \eref{Helmholtz} with \eref{FForm} can be written as
\begin{eqnarray*}
	u(x)&=-\sum_{j=1}^M (\lambda_j+\eta_j\cdot\nabla)\Phi(x; z_j)\\
	&=-\frac{\i}{4}\sum_{j=1}^M\frac{\lambda_j |t_j| H_0^{(1)}(k|t_j|)-k\eta_j\cdot t_j H_1^{(1)}(k|t_j|)}{|t_j|},\quad x\in\mathbb{R}^2\backslash\cup_{j=1}^M\{z_j\},
\end{eqnarray*}
or
\begin{eqnarray*}
	u(x)&=-\sum_{j=1}^M (\lambda_j+\eta_j\cdot\nabla)\Phi(x; z_j)\\
	&= -\frac{1}{4\pi}\sum_{j=1}^M\frac{\e^{\i k|t_j|}}{|t_j|^3}\left(\lambda_j |t_j|^2+\eta_j\cdot t_j\left(\i k|t_j|-1\right)\right),\quad x\in\mathbb{R}^3\backslash\cup_{j=1}^M\{z_j\},
\end{eqnarray*}
where $H^{(1)}_1$ is the Hankel function of the first kind of order one and $t_j:=x-z_j, j=1,\cdots, M$.
Thus, the Dirichlet data on the measurement surface is $u_\Gamma:=u|_\Gamma$.  By straightforward calculations, it can be shown that the Neumann data on the measurement surface are given by
\begin{eqnarray*}
     \partial_{\nu} u_\Gamma(x) &:=\nu(x)\cdot\nabla u(x)\\
     &=\frac{\i k\nu(x)\cdot}{4}\sum_{j=1}^M\frac{1}{|t_j|^3}\Big((\lambda_j t_j+\eta_j) H_1^{(1)}(k|t_j|)|t_j|^2 \\
       &\quad+\eta_j\cdot t_j\bigg(k |t_j| H_0^{(1)}(k|t_j|)-2 H_1^{(1)}(k|t_j|)\bigg) t_j\Big),\quad x\in\Gamma\subset\mathbb{R}^2,
\end{eqnarray*}
or
\begin{eqnarray*}
    \partial_\nu u_\Gamma(x)&:=\nu(x)\cdot \nabla u(x)\\
   &=-\frac{\nu(x)\cdot}{4\pi}\sum_{j=1}^M\frac{\e^{\i k |t_j|}}{|t_j|^5}\Big((\lambda_j t_j+\eta_j)(\i k|t_j|-1)|t_j|^2 \\
    &\quad -\eta_j\cdot t_j(k^2 |t_j|^2+3\i k|t_j|-3)t_j\Big),\quad x\in\Gamma\subset\mathbb{R}^3,
\end{eqnarray*}
where $\nu$ denotes the unit outward normal to $\Gamma$.

In the 2D case, the measurement curve was chosen as the circle centered at the origin with radius of $R$ and $200$ equally spaced measurement angles.  In the 3D case,  the measurement surface is chosen as the sphere centered at the origin with radius $R=6$ and 1806 uniformly distributed observation directions.
%If not elsewhere specified, $N_m=200$ was used in the experiments.

To test the stability of the inversion schemes, in addition to the usual rounding errors, we also contaminate the forward Cauchy data $u_\Gamma$ and $\partial_{\nu}u_\Gamma$ by adding random noise. The noisy Cauchy data $u_\Gamma^\epsilon$ and $\partial_{\nu}u^\epsilon_\Gamma$ were given by the following formulae:
\begin{eqnarray*}
& u^\epsilon_\Gamma:=u_\Gamma+\epsilon r_1|u_\Gamma|\mathrm{e}^{\mathrm{i}\pi r_2},\\
& \partial_\nu u^\epsilon_\Gamma:=\partial_\nu u_\Gamma+\epsilon r_1|\partial_\nu u_\Gamma|\mathrm{e}^{\mathrm{i}\pi r_2}
\end{eqnarray*}
where $r_1$ and $r_2$ are two uniformly distributed random numbers, both ranging from $-1$ to $1$, and $\epsilon>0$ signifies the noise level.

\subsection{Two-dimensional examples}

First we will present some validating examples in two dimensions. In the following figures regarding geometrical settings of the problem, the small red points denote the exact locations of the source points, the blue circle denotes the measurement curve and the black dashed line denotes the boundary of the global sampling domain. In all these examples, the global sampling domain was chosen as the rectangular domain with $100\times 100$ uniformly distributed sampling points. The local sampling domains were chosen as  the square domains centered at the significant global maximizers with side-length  $2\pi/k$  and $40\times 40$ uniformly distributed sampling points.
\\

\begin{figure}
	\centering
	\subfloat[]{\includegraphics[width=0.45\textwidth]{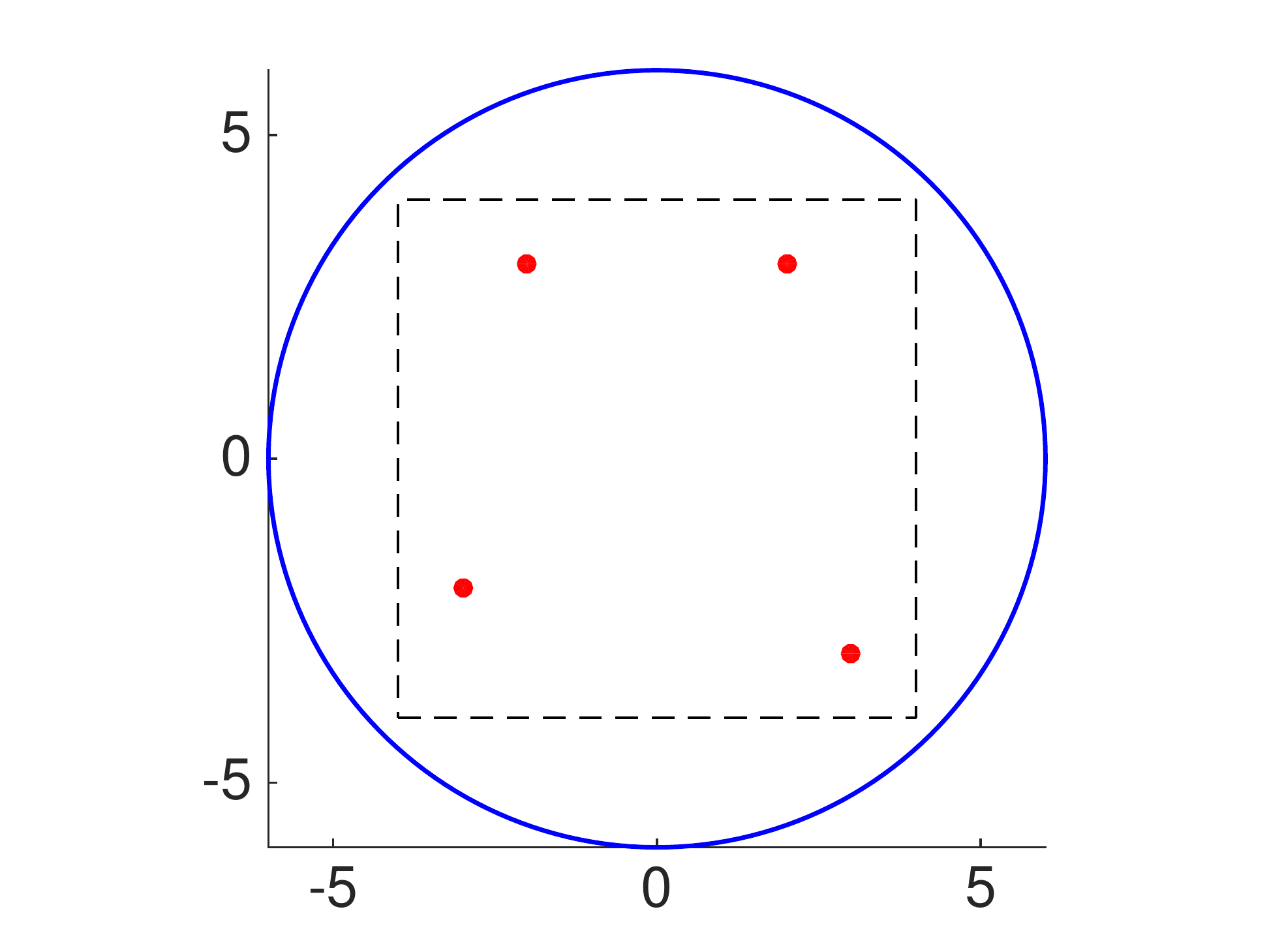}}
	\subfloat[]{\includegraphics[width=0.45\textwidth]{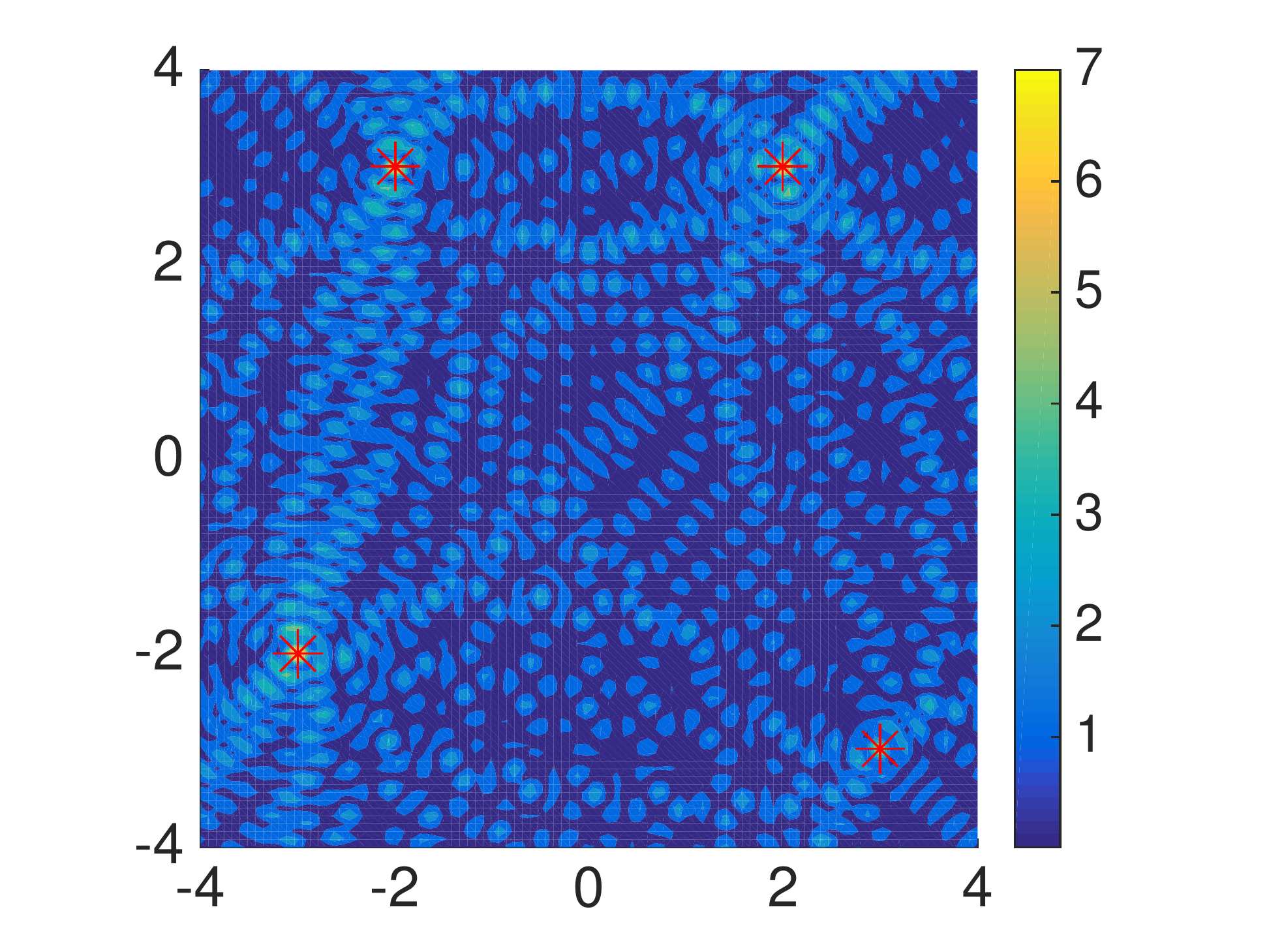}}\\
        \subfloat[]{\includegraphics[width=0.45\textwidth]{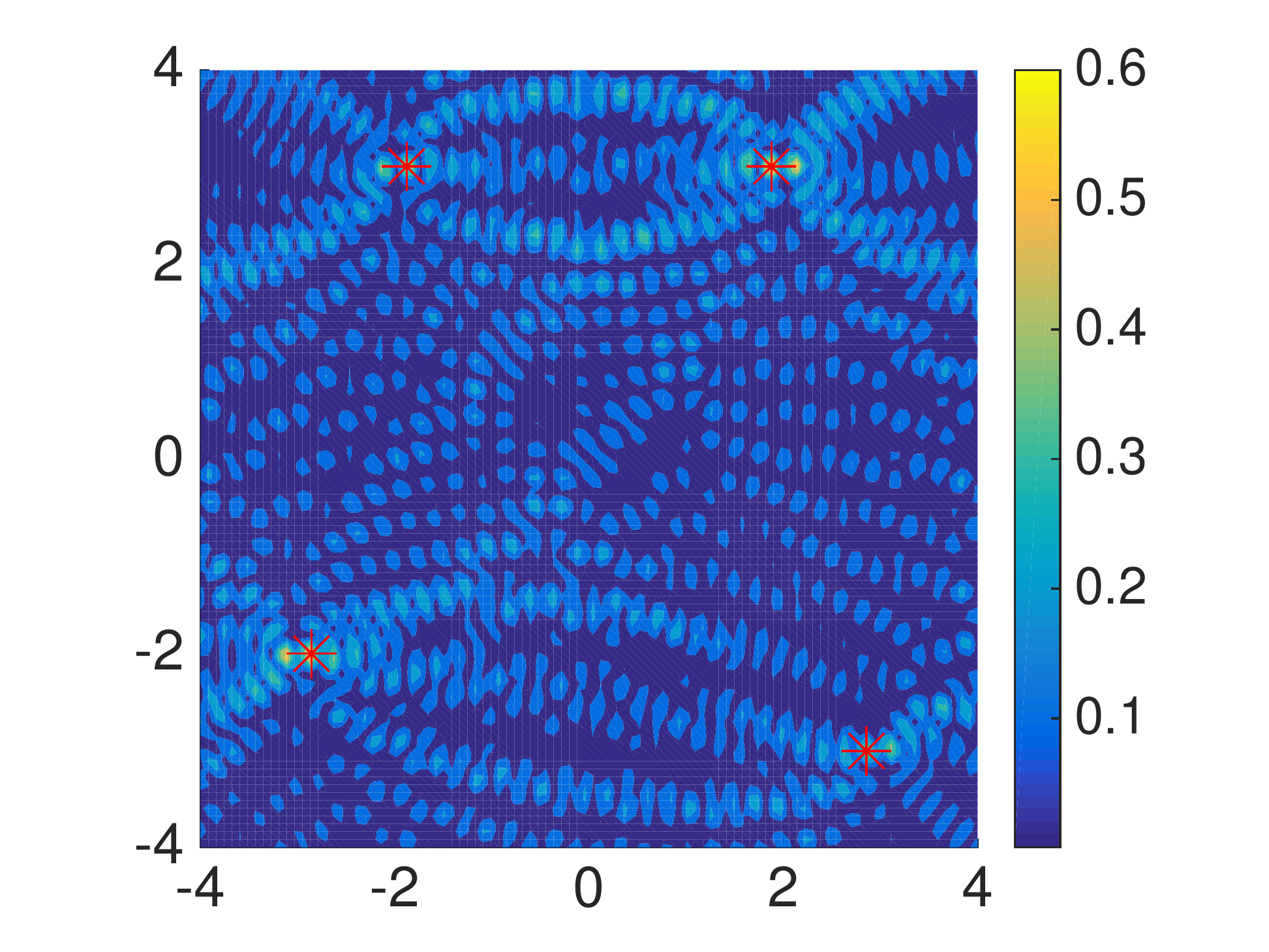}}
        \subfloat[]{\includegraphics[width=0.45\textwidth]{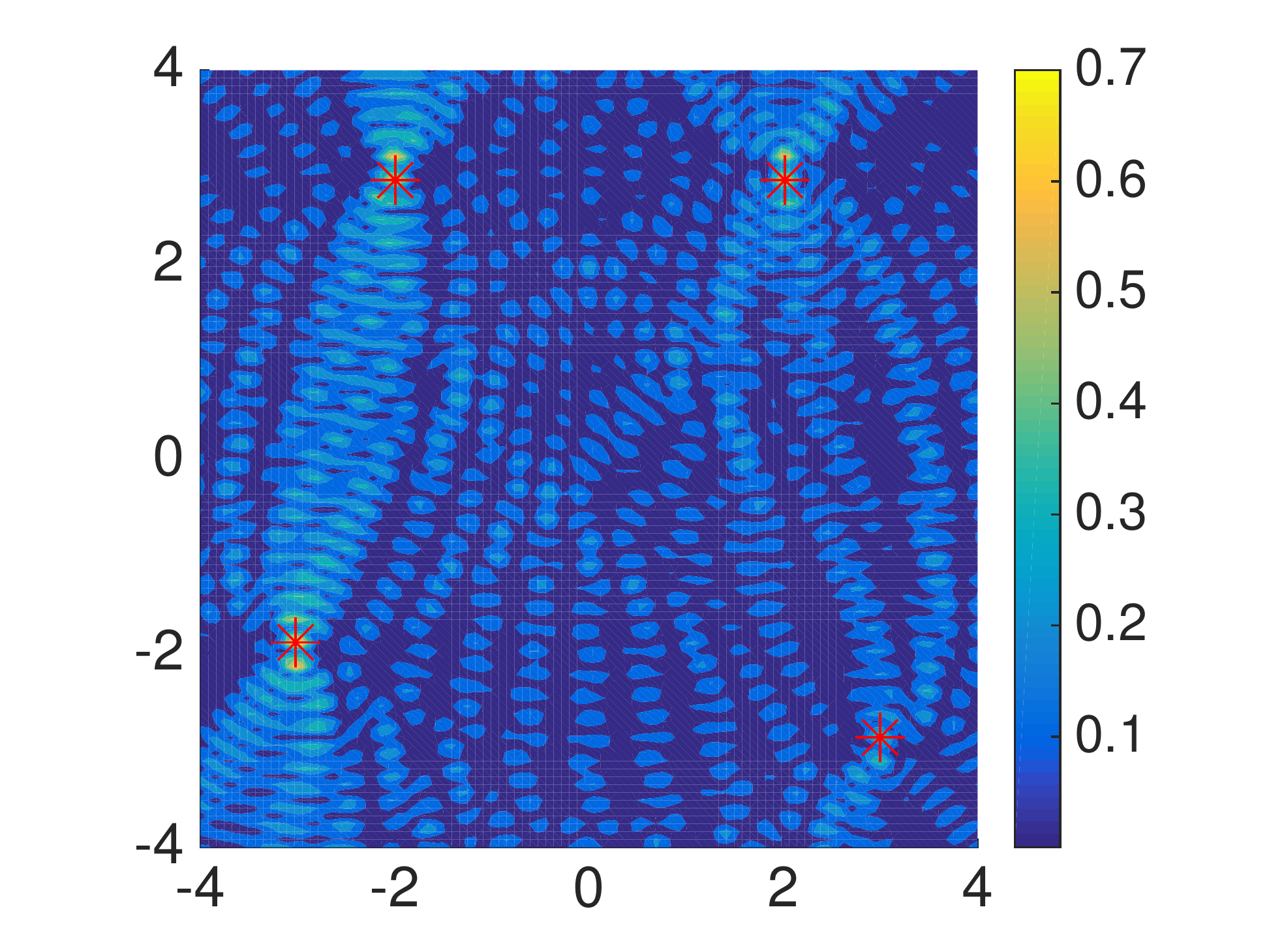}}
	\caption{Geometry setting and images of $z\mapsto |I_{2,\ell}(z)|$ for reconstructing four monopoles in 2D (the points where each indicator attains its largest four extrema are marked by the small red stars). (a) problem geometry (b) $\ell=0$ (c) $\ell=1$ (d) $\ell=2$. }
	\label{fig:example_2D_monopole}
\end{figure}

\noindent {\bf Example 1}\quad In the first example, we will recover the locations of four monopolar source points, see Figure \ref{fig:example_2D_monopole}(a) for the geometry setting of the problem. The wavenumber is chosen as $k=15$, the radius of the measurement circle is $R=6$, the noise level is $\epsilon=5\%$ and the global sampling domain is $[-4, 4]\times[-4,4]$. Some parameters of this example are presented in Table \ref{tab:example_2D_monopole}. Images of indicators $|I_{2,\ell}(z)|, \ell=0,1,2$, are plotted in Figure  \ref{fig:example_2D_monopole}(b), (c) and (d), which show that each of $|I_{2,\ell}(z)|, \ell=0,1,2,$ has four significant local maximizers, which complies well with
the exact locations. \\

\begin{table}
\caption{\label{tab:example_2D_monopole} Reconstruction of monopoles in 2D.}
\begin{indented}\item[]
\begin{tabular}{cccc}
  \br
  & & \multicolumn{2}{c}{\underline{\qquad\qquad Location \qquad\quad}} \\
   Type \& Label & Intensity & Exact & Reconstructed  \\
  \mr
monopole  1 & 9   & $(2, 3)$     &  $(1.9598, 2.9625)$ \\
monopole  2 & 8   & $(-3, -2)$  &  $(-2.9590,-1.9633)$ \\
monopole  3 & 8   & $(-2, 3)$  & $(-1.9529, -2.9495)$ \\
monopole  4 & 7   & $(3, -3)$  & $(2.9495, -2.9495)$ \\
  \br
\end{tabular}
\end{indented}
\end{table}

\begin{figure}
	\centering
	\subfloat[]{\includegraphics[width=0.45\textwidth]{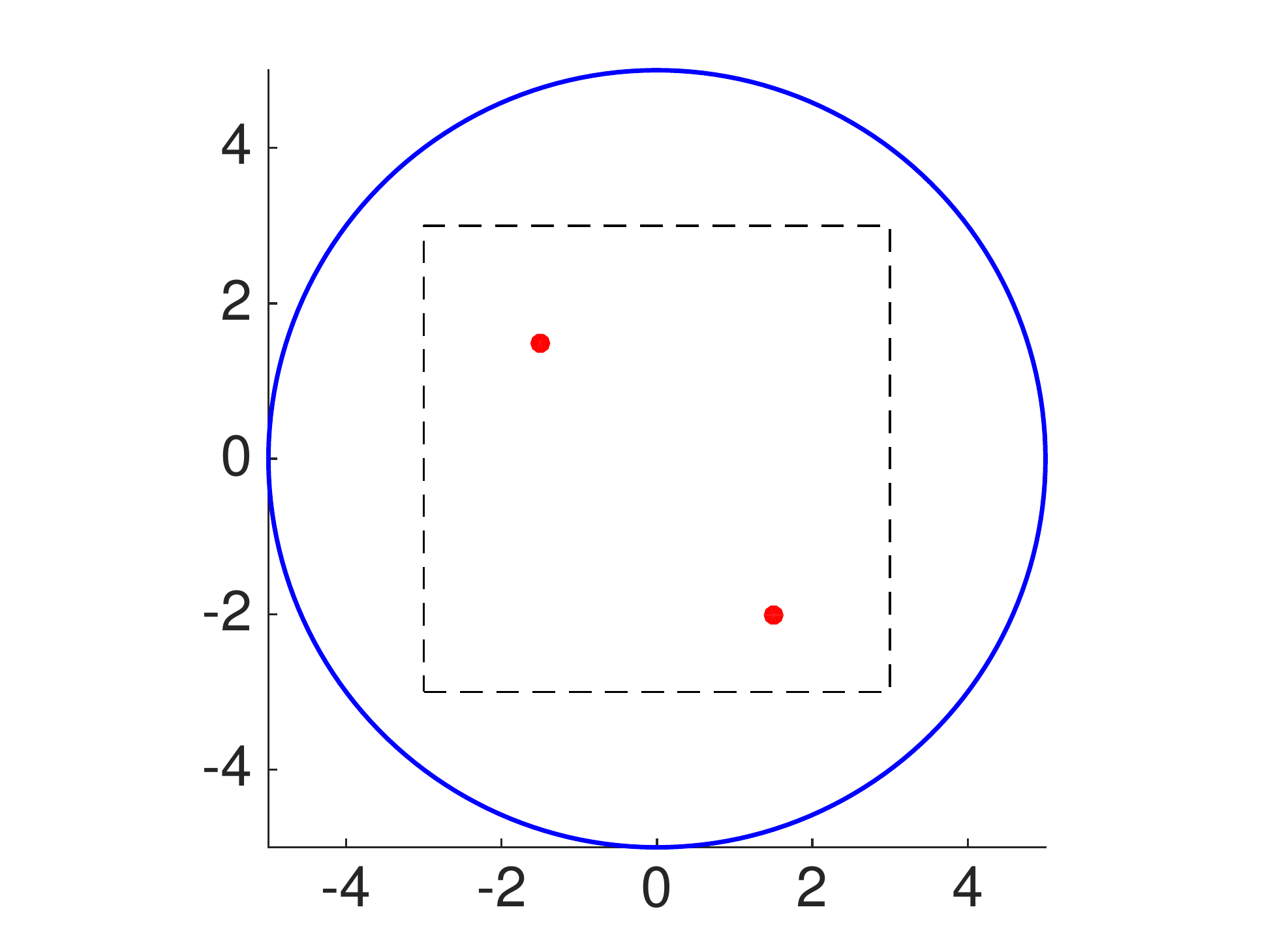}}
	\subfloat[]{\includegraphics[width=0.45\textwidth]{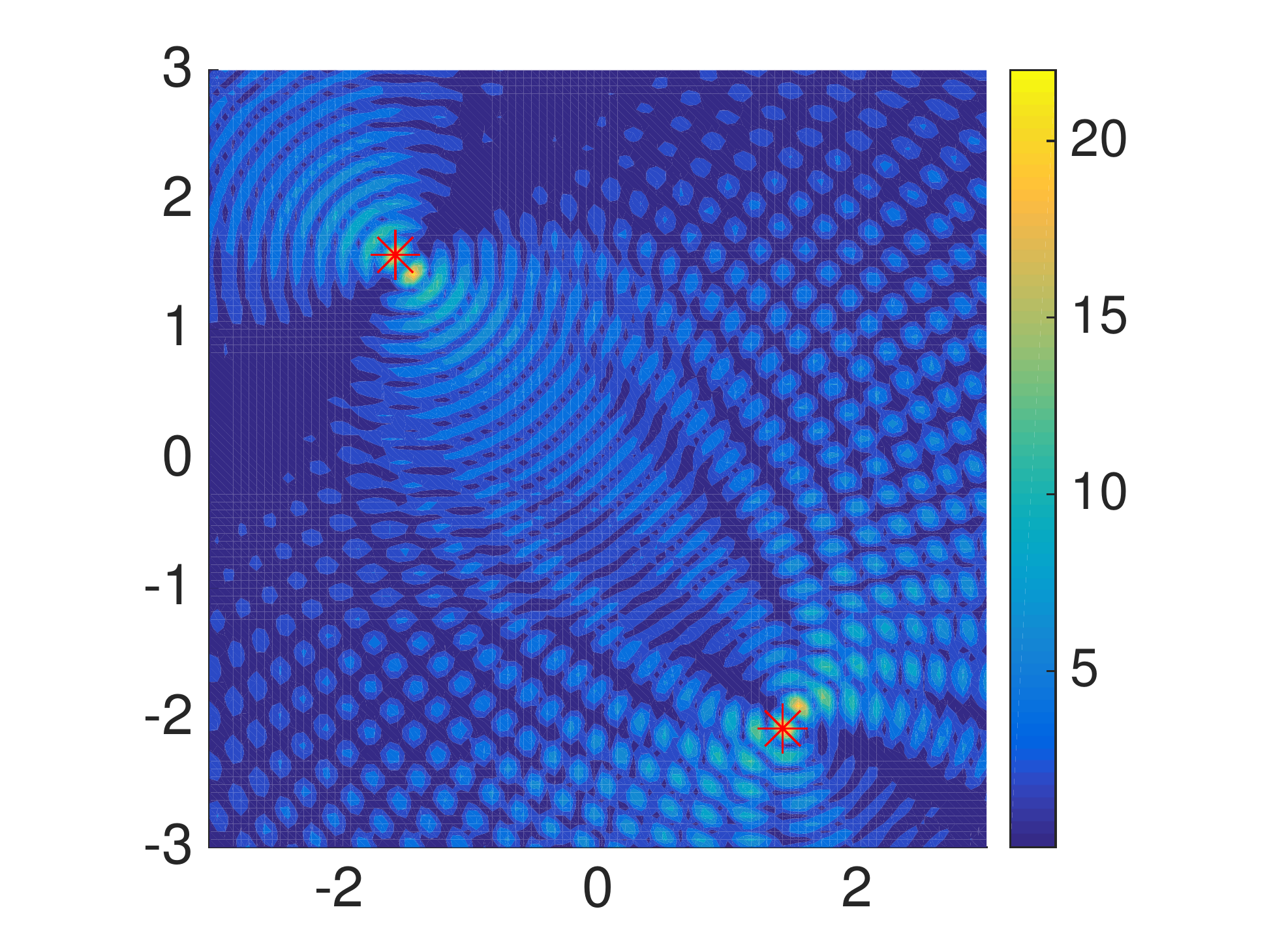}}\\
        \subfloat[]{\includegraphics[width=0.45\textwidth]{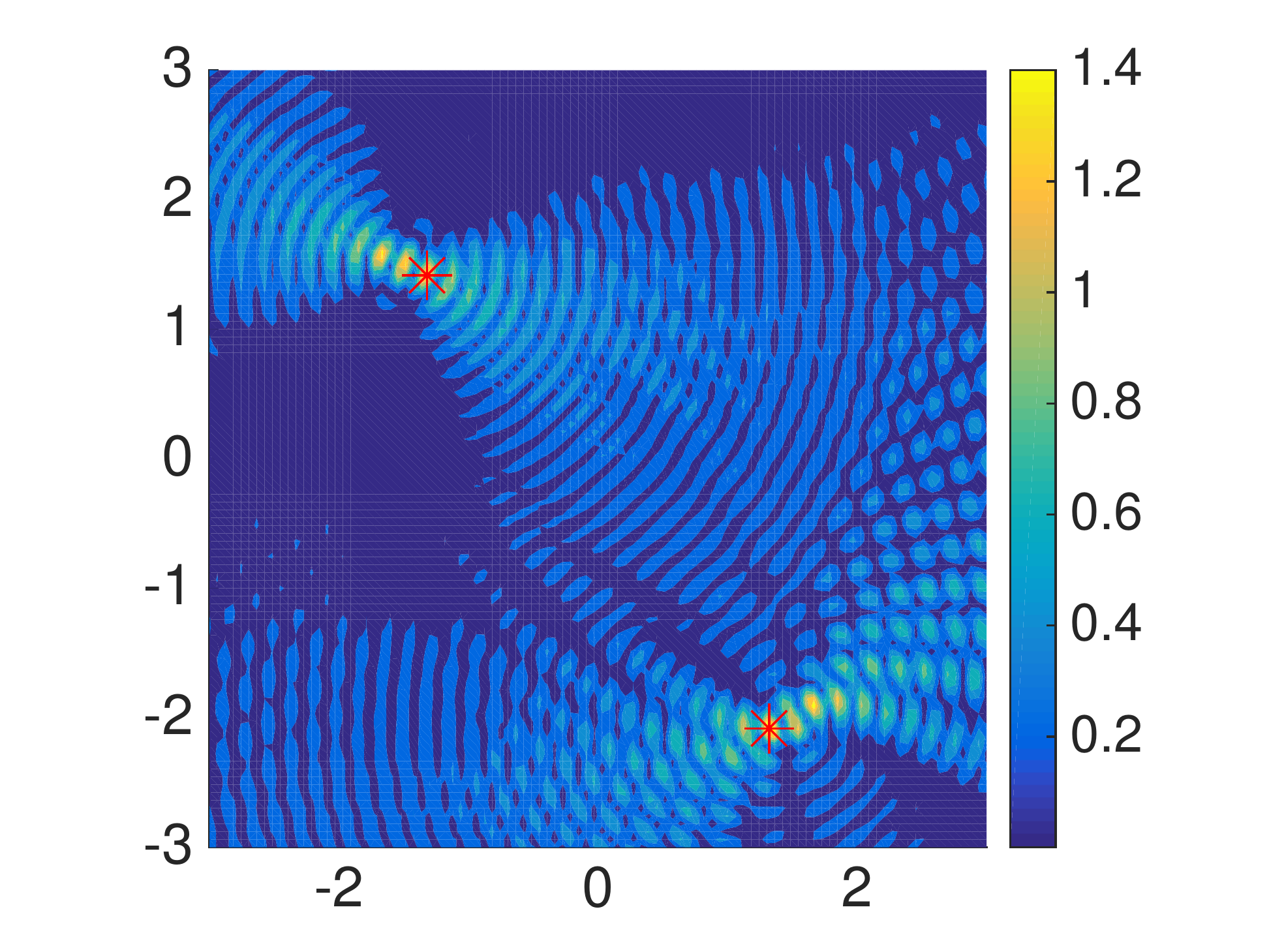}}
        \subfloat[]{\includegraphics[width=0.45\textwidth]{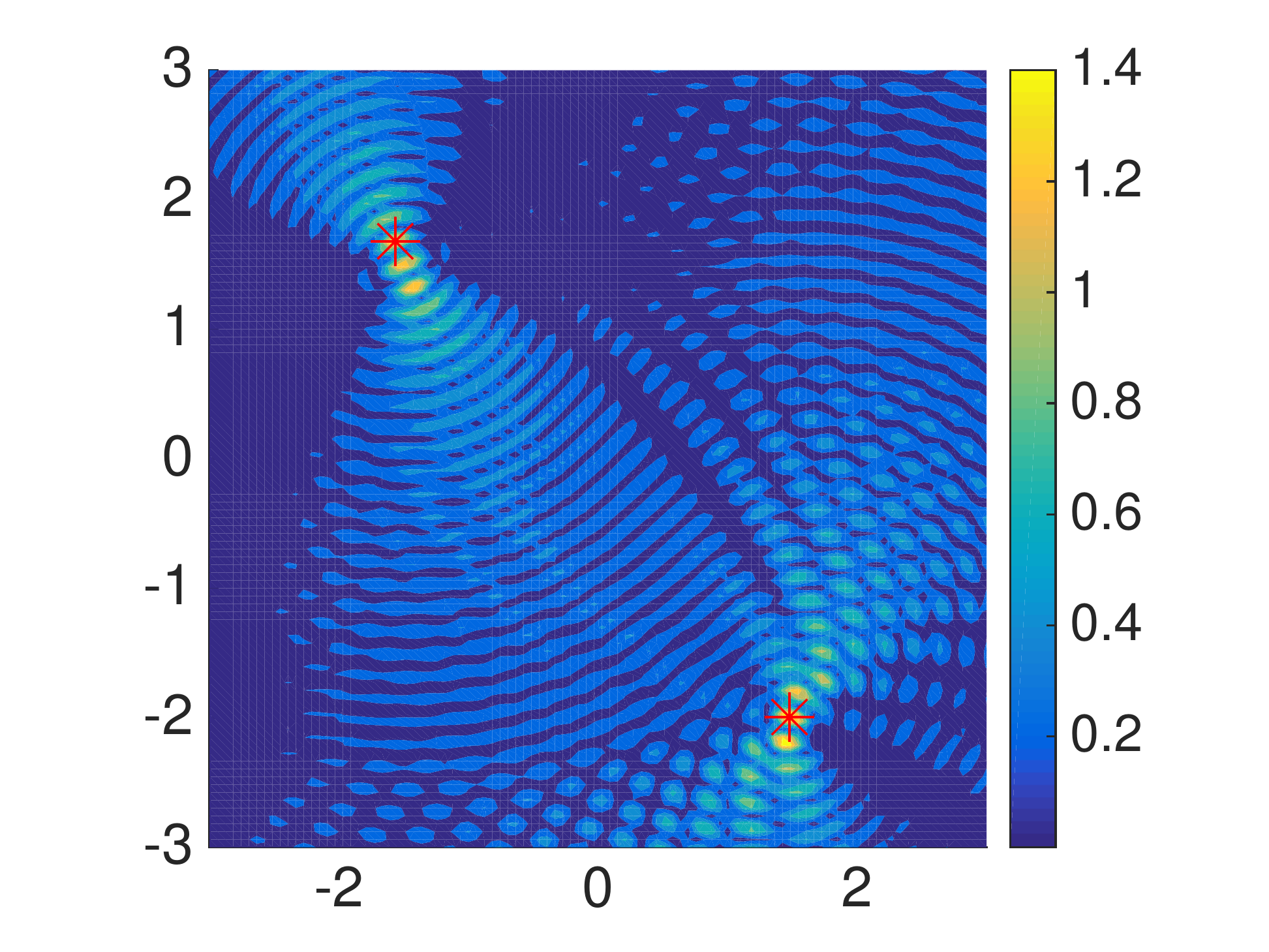}}
	\caption{Geometry setting and images of $z\mapsto |I_{2,\ell}(z)|$ for reconstructing two dipoles in 2D (the points where each indicator attains its largest two extrema are marked by the small red stars). (a) problem geometry (b) $\ell=0$ (c) $\ell=1$ (d) $\ell=2$. }
	\label{fig:example_2D_dipole}
\end{figure}

\noindent {\bf Example 2}\quad Next we will try to reconstruct the locations of two dipole sources, see Figure \ref{fig:example_2D_dipole}(a) for the problem geometry. In this example, $k=18, R=5$ and $\epsilon=5\%$ were used, and the global sampling domain is $[-3, 3]\times[-3,3]$. . Other relevant parameters are listed in Table \ref{tab:example_2D_dipole}. Images of indicators $|I_{2,\ell}(z)|, \ell=0,1,2$, are plotted in Figure  \ref{fig:example_2D_dipole}(b), (c) and (d), which show that each of $|I_{2,\ell}(z)|, \ell=0,1,2,$ has two significant local maximizers, which are consistent with exact locations. It is worth remarking that the dipole sources will incur several spurious spikes nearby, among which we always pick up the point with largest indicator value and rule out the others due to the assumption of sparse source distribution.
\\

\begin{figure}
	\centering
	\subfloat[]{\includegraphics[width=0.45\textwidth]{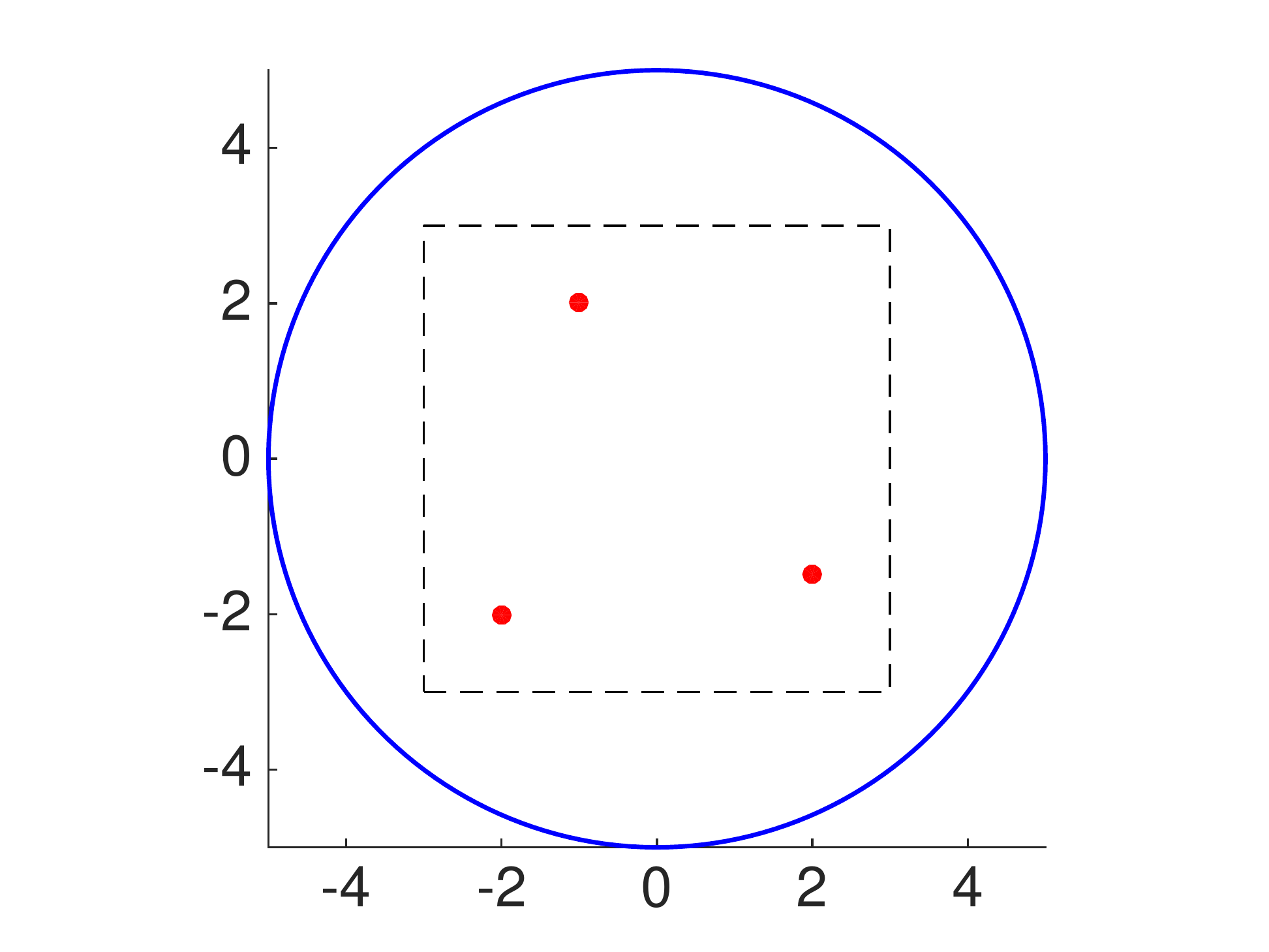}}
	\subfloat[]{\includegraphics[width=0.45\textwidth]{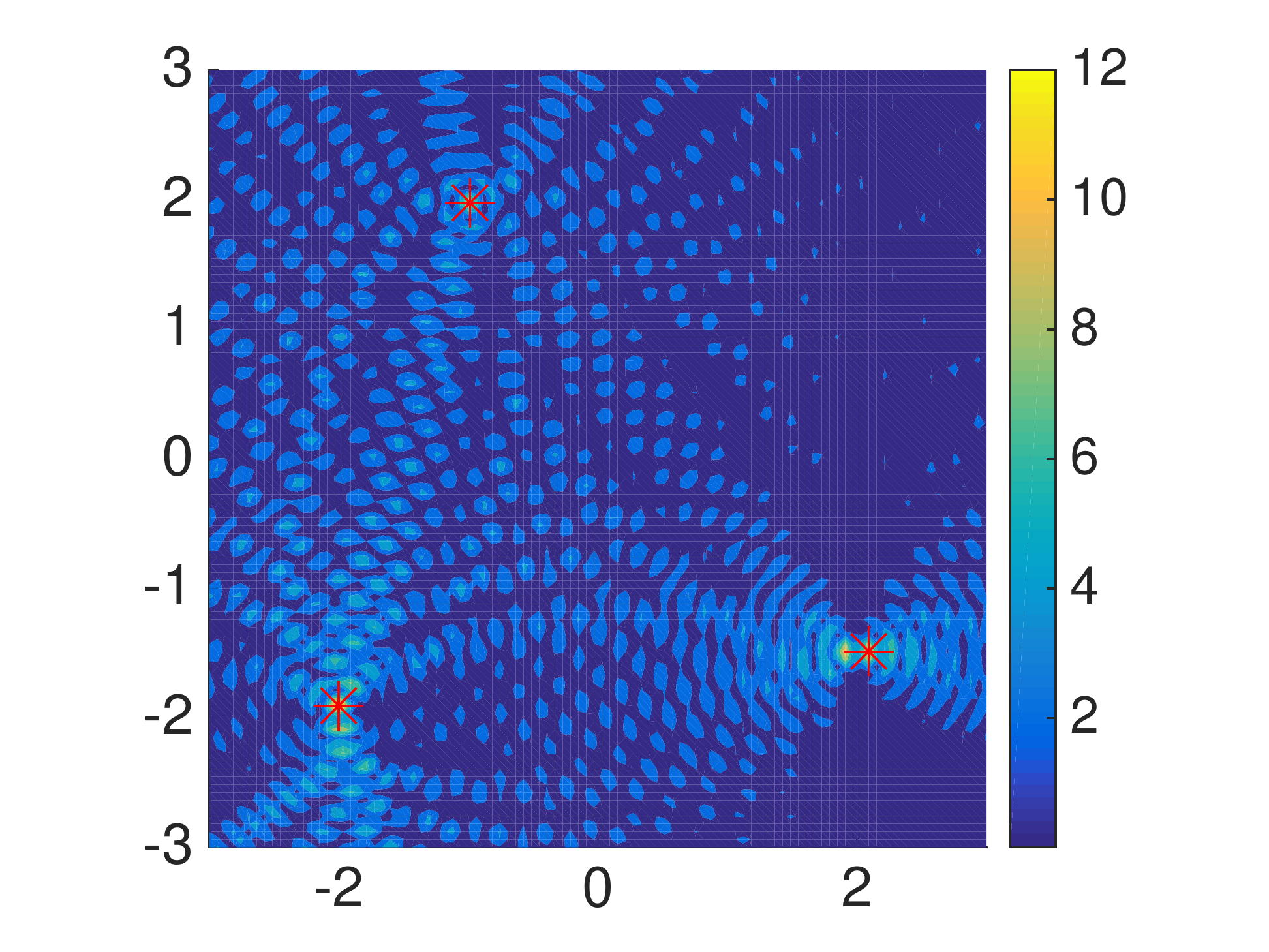}}\\
        \subfloat[]{\includegraphics[width=0.45\textwidth]{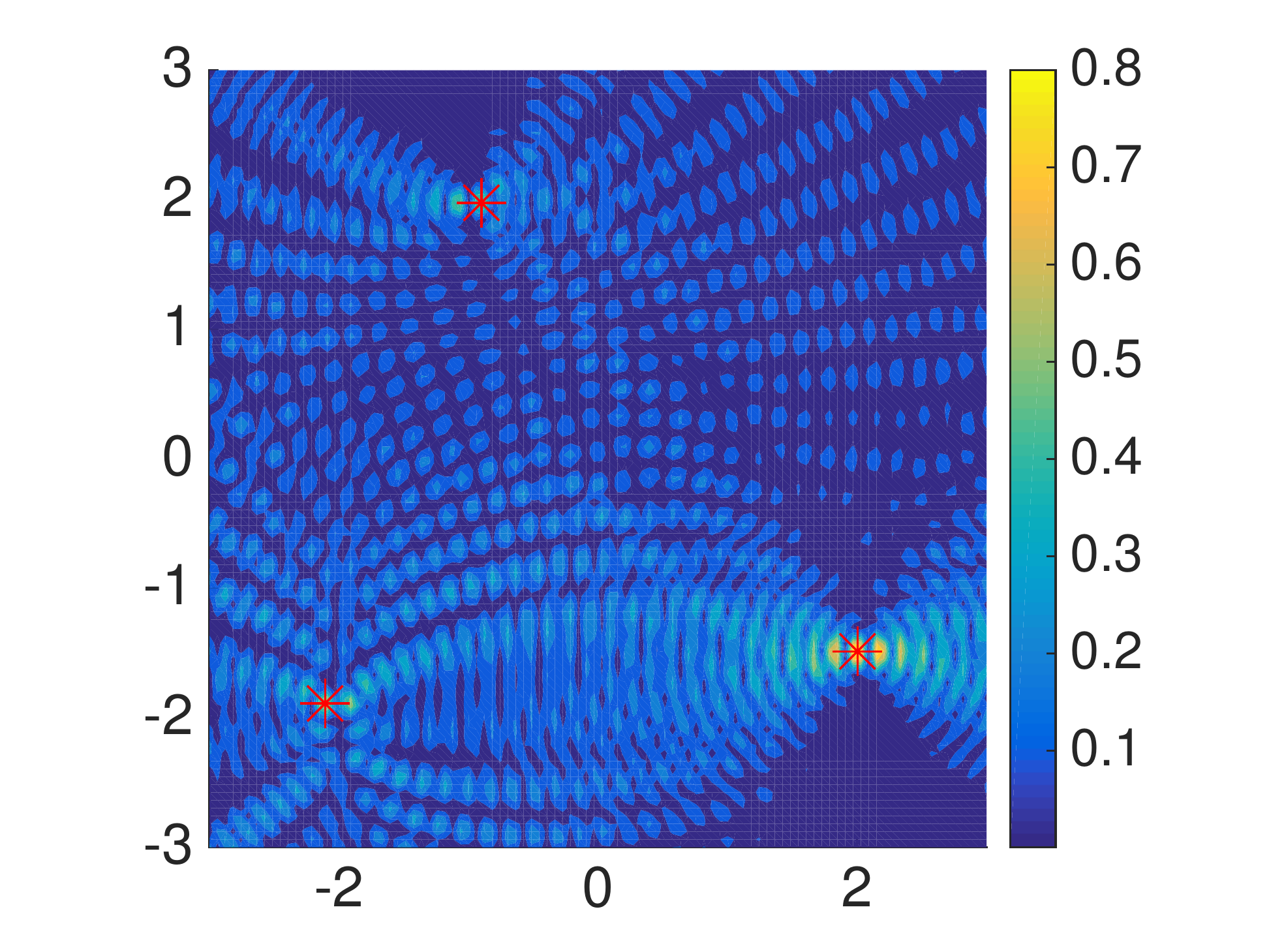}}
        \subfloat[]{\includegraphics[width=0.45\textwidth]{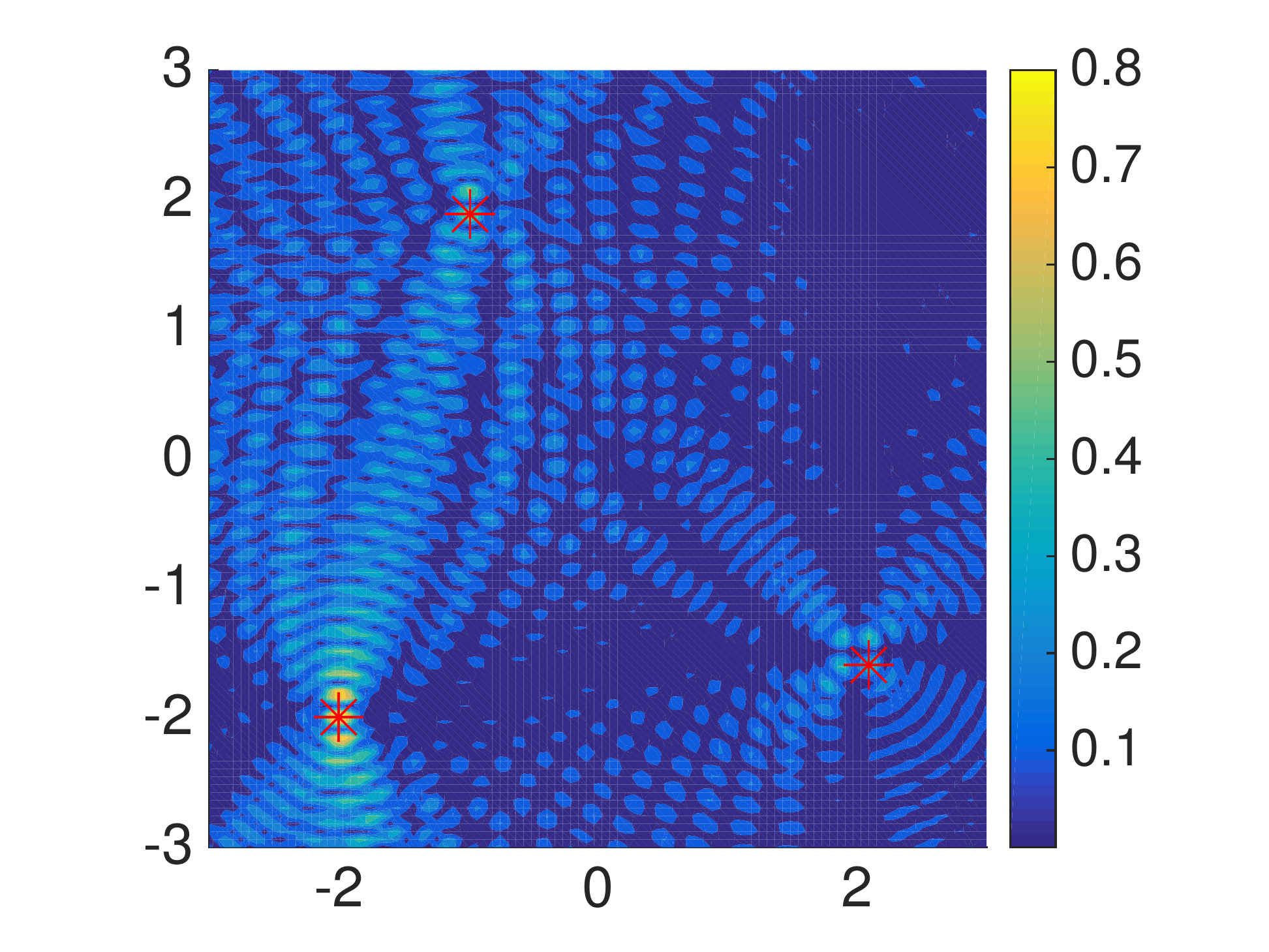}}
	\caption{Geometry setting and images of $z\mapsto |I_{2,\ell}(z)|$ for reconstructing a monopole and two dipoles in 2D (the points where each indicator attains its largest three extrema are marked by the small red stars). (a) problem geometry (b) $\ell=0$ (c) $\ell=1$ (d) $\ell=2$. }
	\label{fig:example_2D_multipole}
\end{figure}

\noindent {\bf Example 3}\quad In this example, we aim to reconstruct a source consisting of a monopolar component and two dipolar components. The geometry configuration of this problem is shown in Figure \ref{fig:example_2D_multipole}(a). The exact locations and intensities of the point sources are listed in Table \ref{tab:example_2D_multipole}. Here we use $k=20, R=5, \epsilon=5\%$ and the global sampling domain $[-3, 3]\times[-3,3]$. . Images of indicators $|I_{2,\ell}(z)|, \ell=0,1,2$, are plotted in Figure  \ref{fig:example_2D_multipole}(b), (c) and (d), which show that each of $|I_{2,\ell}(z)|, \ell=0,1,2,$ has three significant local maximizers. The coordinates of these local maximizers are listed in Table \ref{tab:max_points}.

%\begin{table}
%\caption{\label{tab:example_2D_monopole} Some parameters in the reconstruction of two monopoles.}
%\begin{indented}\item[]
%\begin{tabular}{ccccc}
%  \br
% &  & & \multicolumn{2}{c}{\underline{\qquad\qquad Location \qquad\quad}} \\
%   Source  & Type & Intensity & Exact & Reconstructed  \\
%  \mr
%1 &	monopole  & 20   & $(1, 1)$       &  $(0.9666, 0.9631)$ \\
%2 &	monopole  & 22   & $(-1.5, -1)$  &  $(-1.4689,-0.9631)$ \\
%  \br
%\end{tabular}
%\end{indented}
%\end{table}

\begin{table}
\caption{\label{tab:example_2D_dipole} Reconstruction of dipoles in 2D.}
\begin{indented}\item[]
\begin{tabular}{cccc}
  \br
& & \multicolumn{2}{c}{\underline{\qquad\qquad Location \qquad\quad}} \\
 Type \& Label & Intensity & Exact & Reconstructed  \\
  \mr
dipole 1 & $(-\sqrt{2}, \sqrt{2})$  & $(-1.5, -1.5)$ &  $(-1.4386, 1.5109)$ \\
dipole 2 & $(\sqrt{2}, \sqrt{2})$   & $(1.5, -2)$     &  $(1.4157, -2.0534)$ \\
  \br
\end{tabular}
\end{indented}
\end{table}

\begin{table}
\caption{\label{tab:example_2D_multipole} Reconstruction of multipoles in 2D.}
\begin{indented}\item[]
\begin{tabular}{cccc}
  \br
  & & \multicolumn{2}{c}{\underline{\qquad\qquad Location \qquad\quad}} \\
Type \& Label & Intensity & Exact & Reconstructed  \\
  \mr
monopole  & 10       & $(-1, 2)$      &  $(-0.9622, 1.9495)$ \\
dipole 1   & $(1,0)$   & $(2, -1.5)$    &  $(2.0636, -1.5279)$ \\
dipole 2  & $(0,1)$   & $(-2, -2)$     &  $(-2.0374, -1.9293)$\\
  \br
\end{tabular}
\end{indented}
\end{table}

\begin{table}
\caption{\label{tab:max_points} The first three global maximizers of the indicators in Example 3.}
\begin{indented}\item[]
\begin{tabular}{ccc}
			\br
			$|I_{2,0}(z)|$                 & $|I_{2,1}(z)|$                & $|I_{2,2}(z)|$  \\
			\mr
			$(-2.0067,-1.9012)$ & $(-2.1066,-1.8934)$ & $(-1.9989,-1.9933)$   \\
			$(2.0909,-1.4927)$  & $(2.0011,-1.4927)$   & $(2.0988,-1.5982)$    \\
			$(-0.9921,-1.9776)$ & $(-0.9023,1.9854)$  & $(-0.9921,1.8855)$    \\
			\br
\end{tabular}
\end{indented}
\end{table}

%%%%%%%%%%%%%%%%%%%%%%%%%%%%%%%%%%%%%%%%%%%%%%%

\subsection{Three-dimensional examples}

Now we are going to present some 3D examples. In the following 3D examples, the wavenumber is $k=10$ and the global sampling region is chosen as $[-3, 3]^3\subset\mathbb{R}^3$.  For DSM, the sampling region $[-3, 3]^3$ is probed by a uniform grid $G$ of $60\times 60\times 60$ sampling points. For DSM2, the global sampling region $[-3, 3]^3$ is probed by a uniform $30\times 30\times 30$ grid of  sampling points and each local sampling region is a cube with side-length $2\pi/k$ and uniform $20\times 20\times 20$ sampling points. To facilitate the 3D visualizations, some 2D projections (shadows with gray color) are also added to the figures. \\

\begin{figure}
	\centering
	\subfloat[]{\includegraphics[width=0.45\textwidth]{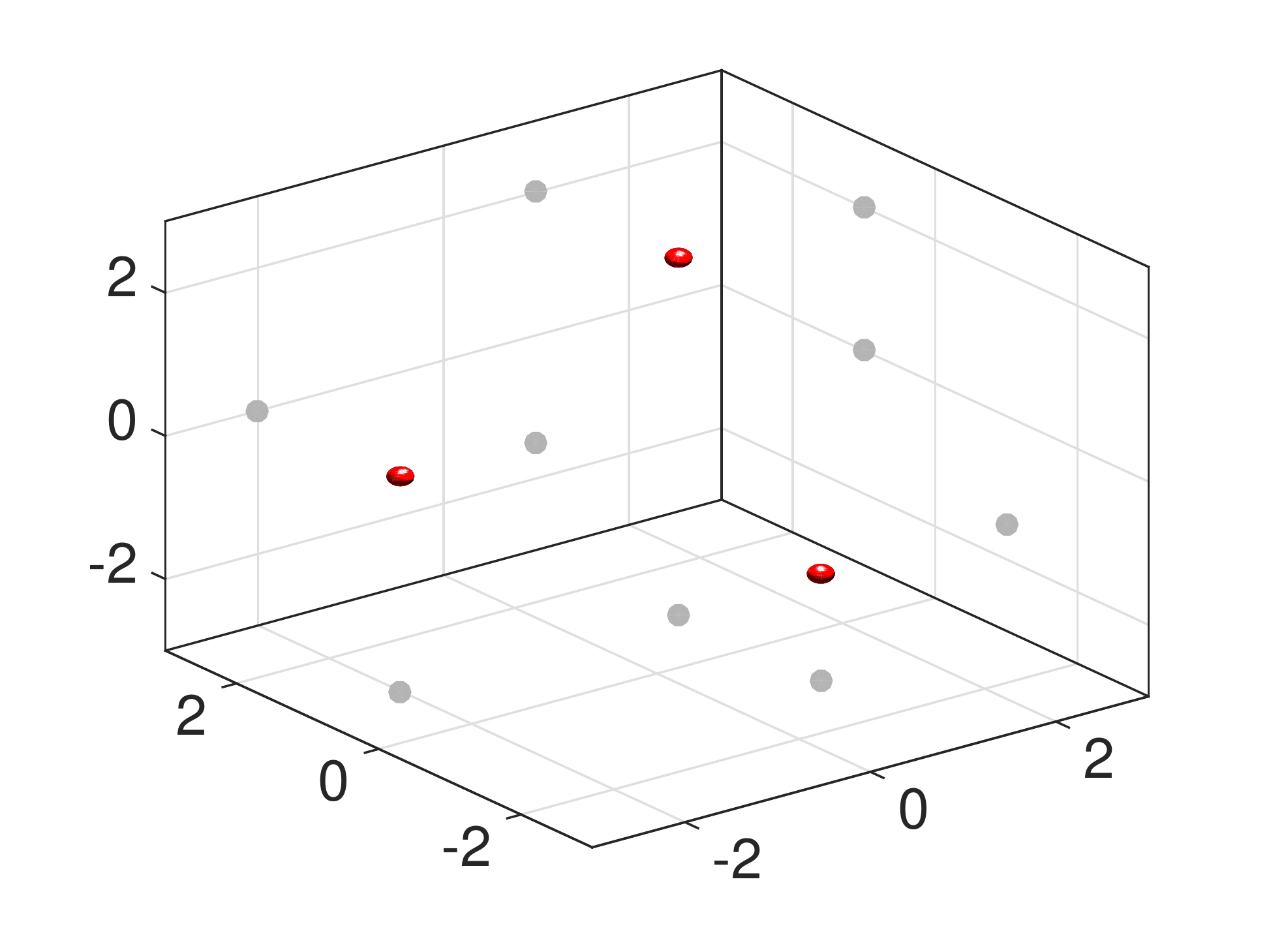}}
	\subfloat[]{\includegraphics[width=0.45\textwidth]{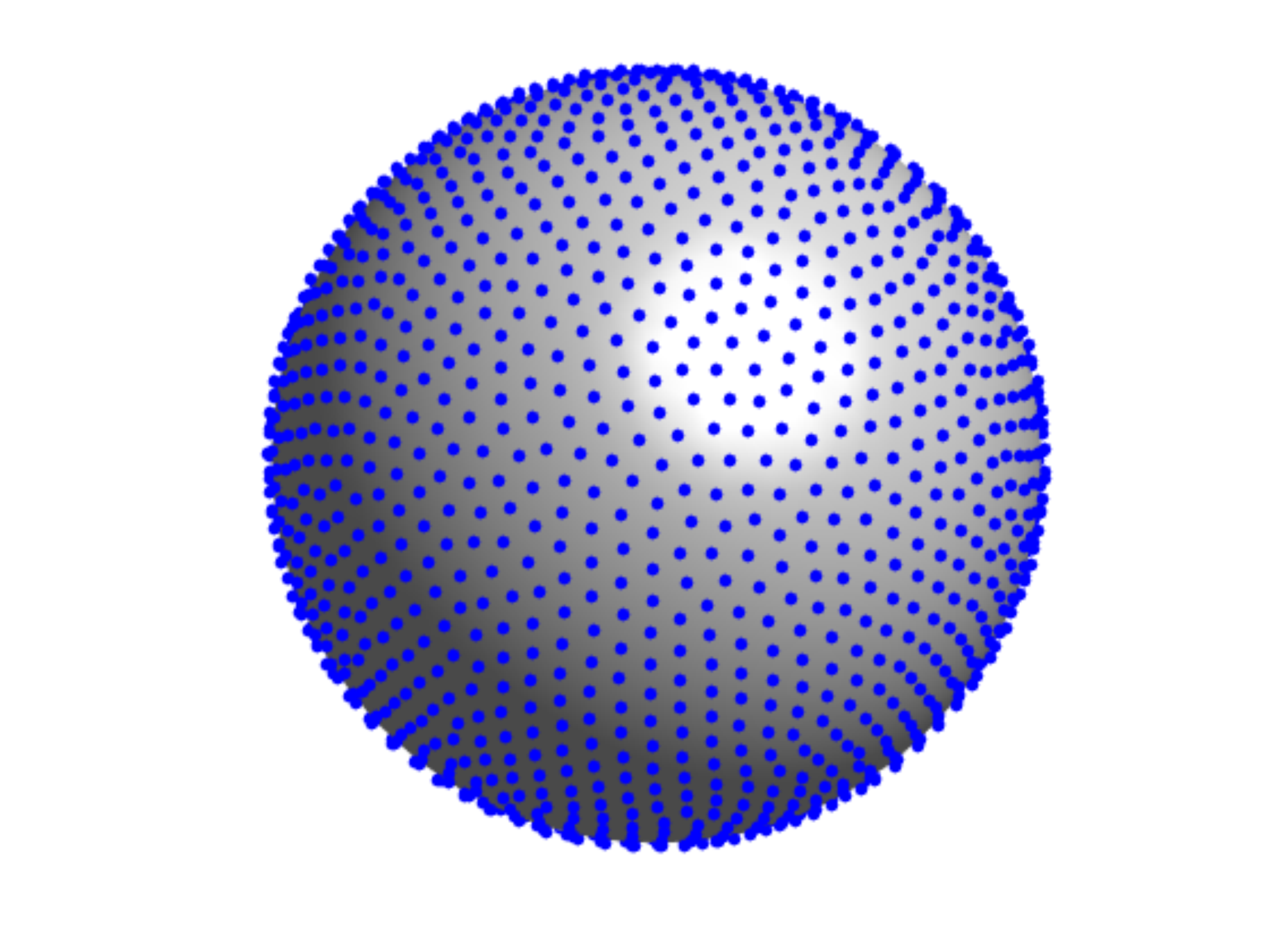}}
	\caption{Geometry setting of the direct problem in 3D. (a) problem geometry (the small red balls denote the exact locations of the source points). (b) configuration of receivers (the 1806 receivers are uniformly distributed on the sphere with radius 6 and denoted by the small blue points).}
	\label{fig:example_3D_monopole_geometry}
\end{figure}

\begin{figure}
	\centering
	\subfloat[]{\includegraphics[width=0.45\textwidth]{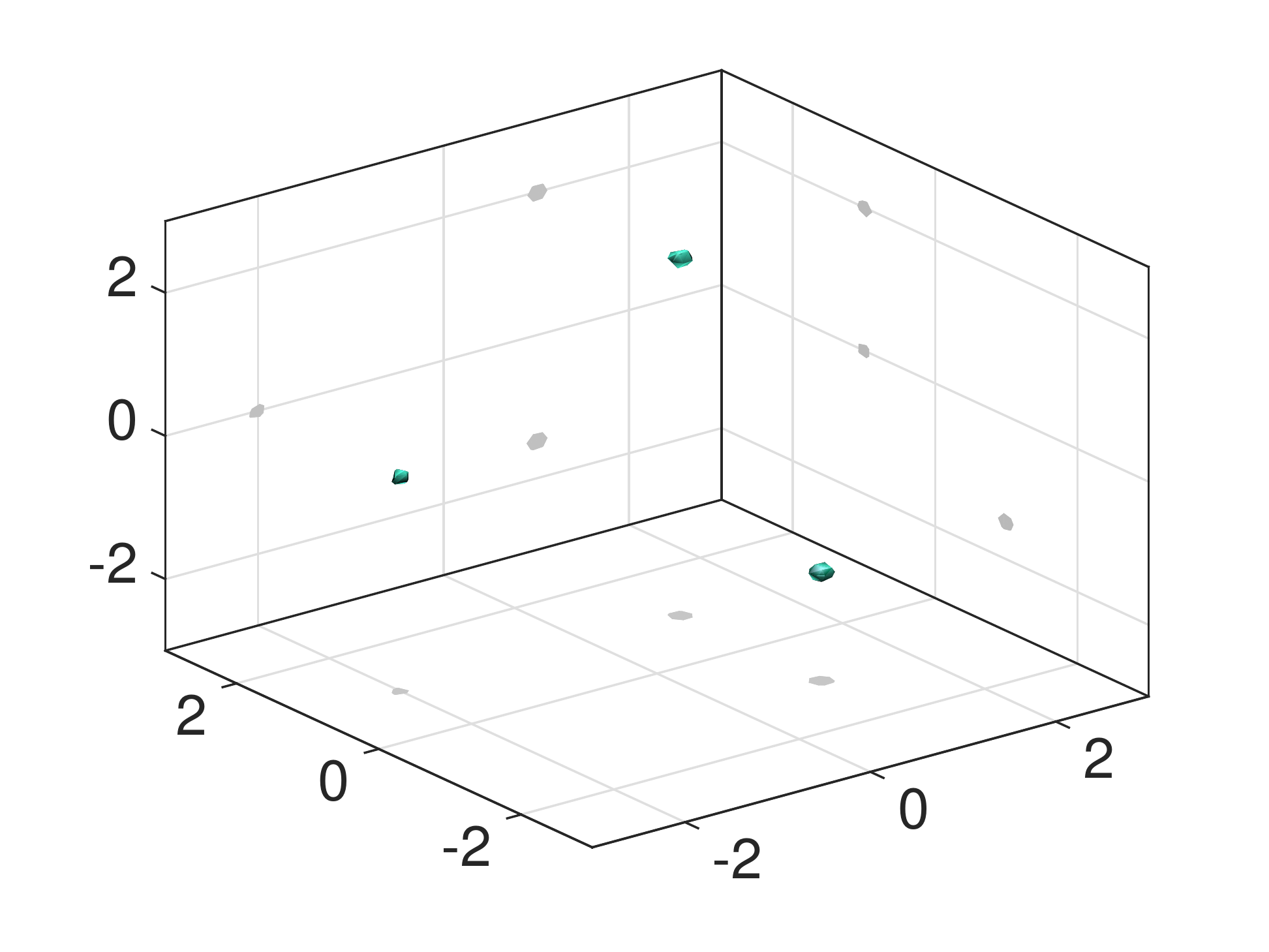}}
	\subfloat[]{\includegraphics[width=0.45\textwidth]{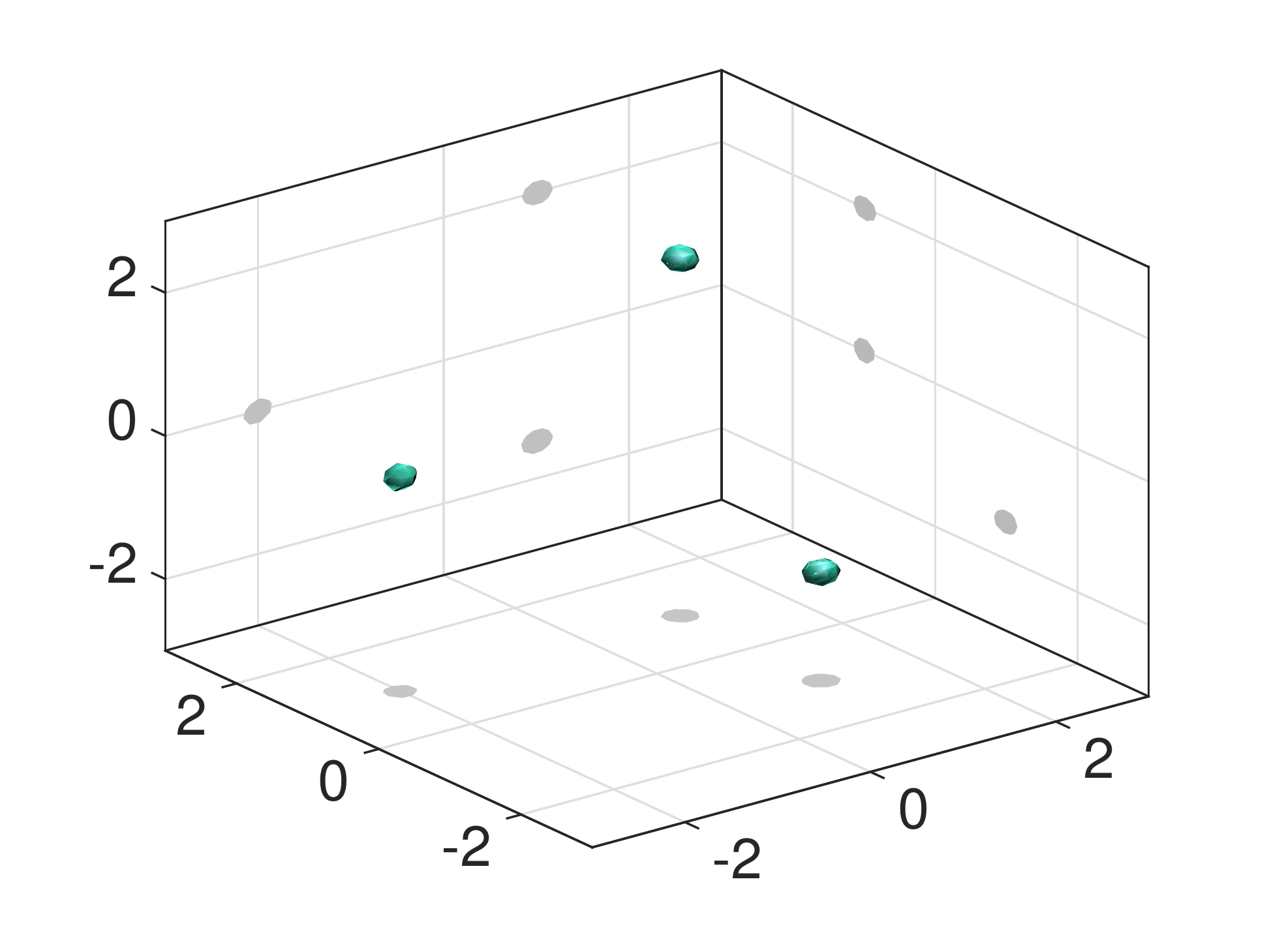}}\\
	\subfloat[]{\includegraphics[width=0.45\textwidth]{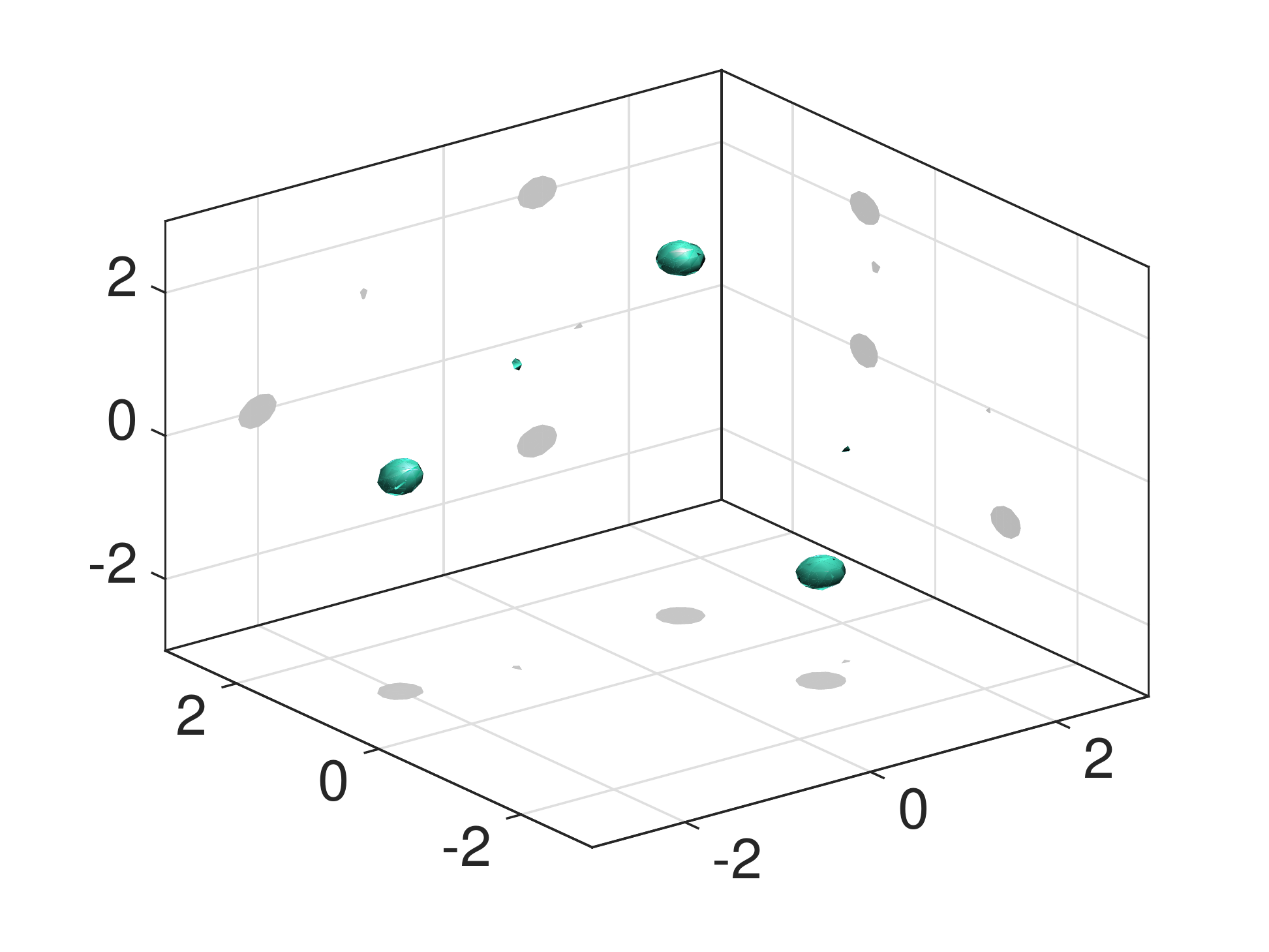}}
	\subfloat[]{\includegraphics[width=0.45\textwidth]{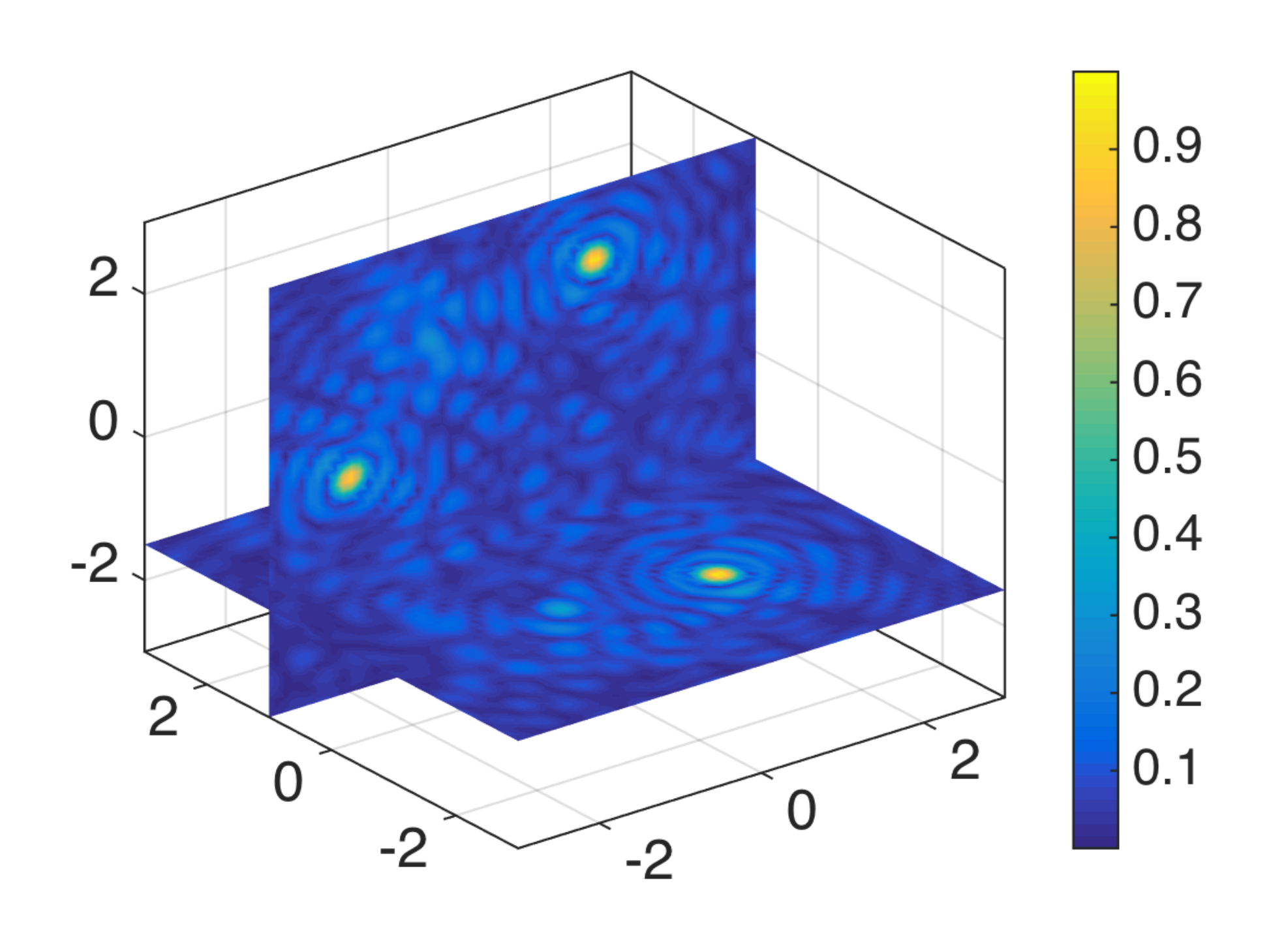}}
	\caption{Iso-surface and cross-section visualizations of  $z\mapsto |I_{3,0}(z)|/\max_{z\in G}|I_{3,0}(z)|$ for reconstructing three monopoles in 3D. The iso-surfaces are drawn with cut-off values $C$. (a)  $C=0.8$  (b)  $C=0.6$  (c)  $C=0.4$ (d) cross-section with slice positions $x_2=1$ and $x_3=-1.5$.}
	\label{fig:example_3D_monopole_reconstruction}
\end{figure}

\noindent {\bf Example 4}\quad First, we want to compare the computational costs of the DSM and DSM2 by reconstructing three monopoles, see Figure \ref{fig:example_3D_monopole_geometry} for the geometrical setup of the forward problem. In this example, the source intensities are $\lambda_1=\lambda_2=\lambda_3=5$ and the noise level is $\epsilon=10\%$.
The reconstructions of DSM are shown in Figure \ref{fig:example_3D_monopole_reconstruction} via plotting the results both in the iso-surface and cross-section formats. For simplicity, we only plot the normalized indicator of $I_{3,0}(z)$. We list the reconstructed locations and the computational CPU time for the reconstructions in Table \ref{tab:example_3D_monopole}. All of the codes in our experiments are written in MATLAB and run on an Intel Core 2.6GHz laptop. As shown in Table \ref{tab:example_3D_monopole}, in comparison to the single-level version DSM, the  ``self-adapted'' two-level version DSM2 could produce more accurate (or at least comparable) result with significantly less computational cost.  \\

\begin{table}
\caption{\label{tab:example_3D_monopole} A comparison of DSM and DSM2 in the reconstruction of three monopoles. Here $T$ denotes the CPU time of computation (in seconds).}
\begin{indented}\item[]
\begin{tabular}{cccc}
  \br
  & & \multicolumn{2}{c}{\underline{\qquad\qquad\qquad Reconstruction \qquad\qquad\qquad}} \\
Label  & Exact location & DSM ($T=576.5$s) & DSM2 ($T=6.8$s)  \\
  \mr
 1 &  $(1, 1, 2)$        & $(0.9661, 0.9661, 1.9831)$     &  $(1.0253, 0.9939, 1.9969)$ \\
 2 &  $(1, -1, -1.5)$  & $(0.9661, -0.9661,-1.4746)$    &  $(0.9939, -0.9939, -1.4889)$ \\
 3 &  $(-2, 1, 0)$      & $(-1.9831, 0.9661,-0.0508)$     &  $(-1.9969, 0.9939, 0.0092)$\\
  \br
\end{tabular}
\end{indented}
\end{table}

\begin{figure}
	\centering
	\subfloat[]{\includegraphics[width=0.45\textwidth]{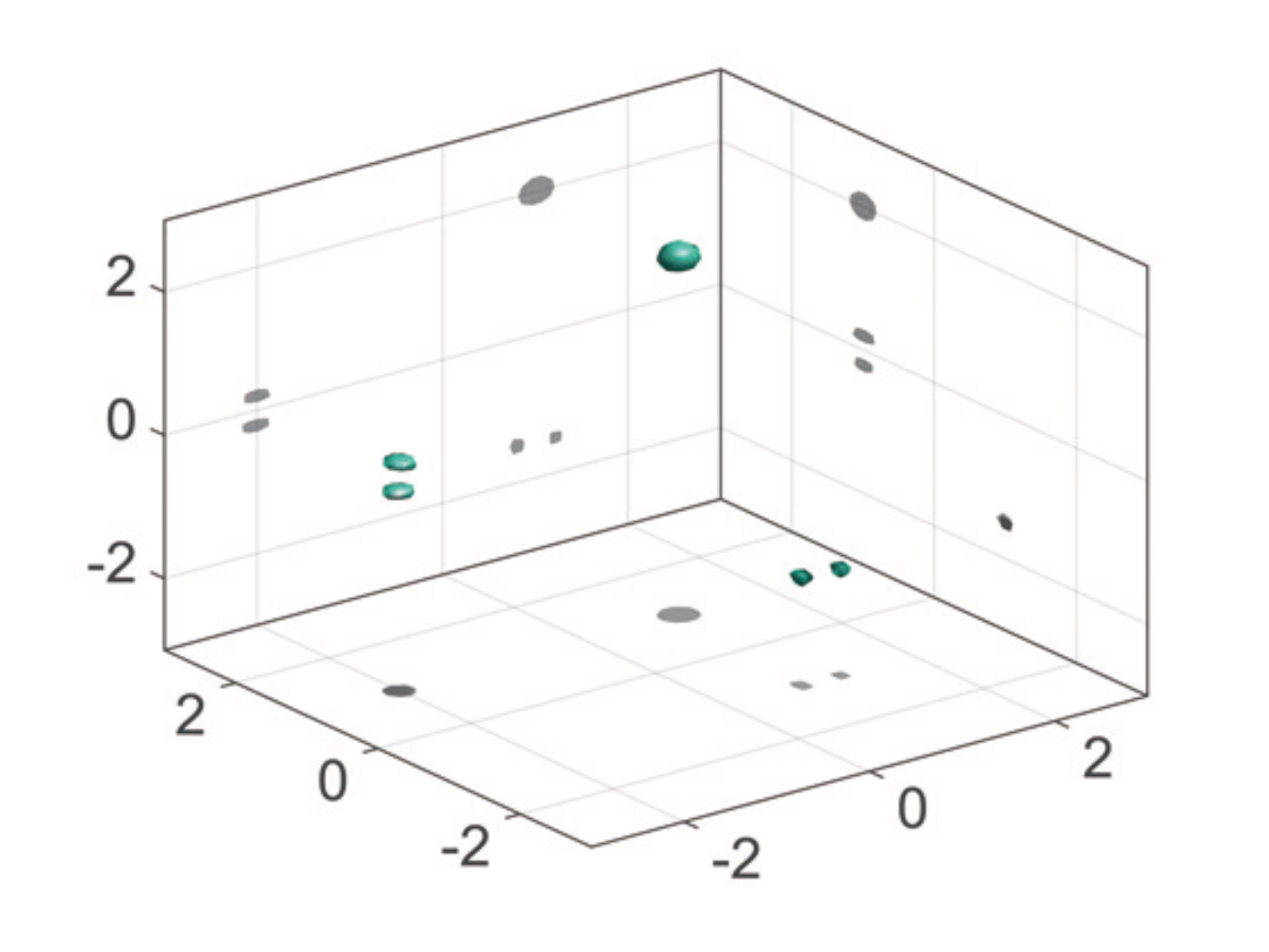}}
	\subfloat[]{\includegraphics[width=0.45\textwidth]{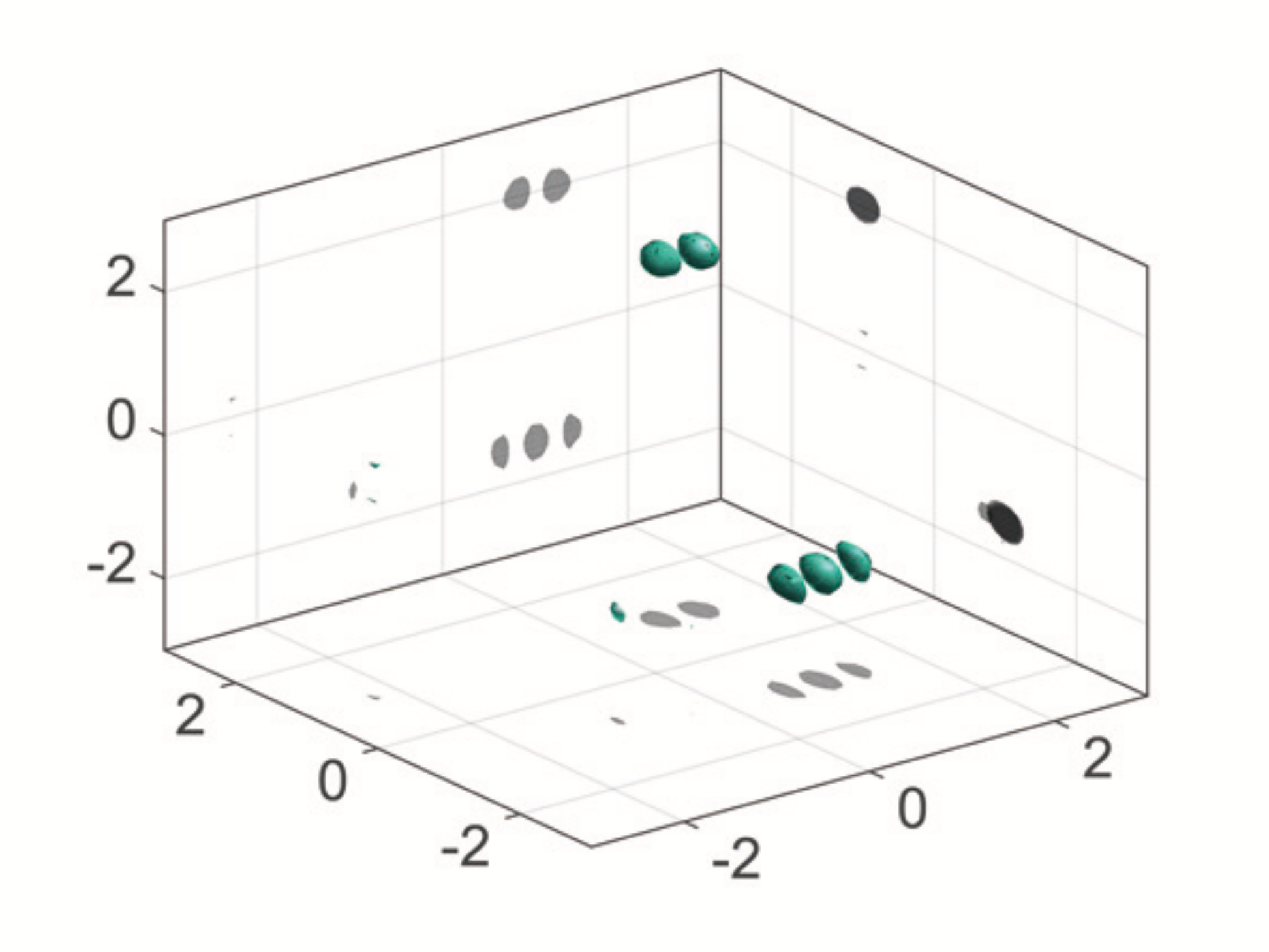}}\\
	\subfloat[]{\includegraphics[width=0.45\textwidth]{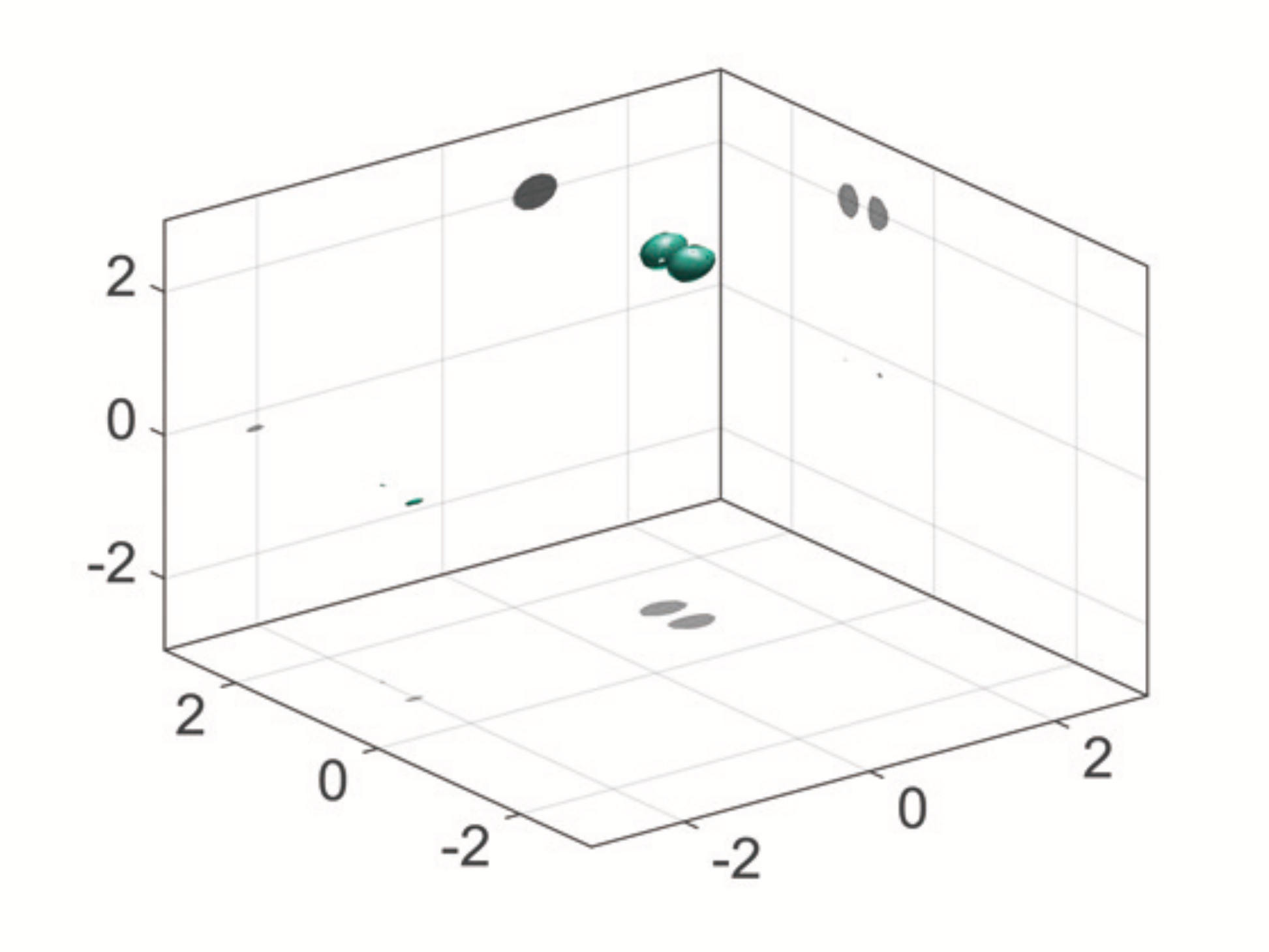}}
	\subfloat[]{\includegraphics[width=0.45\textwidth]{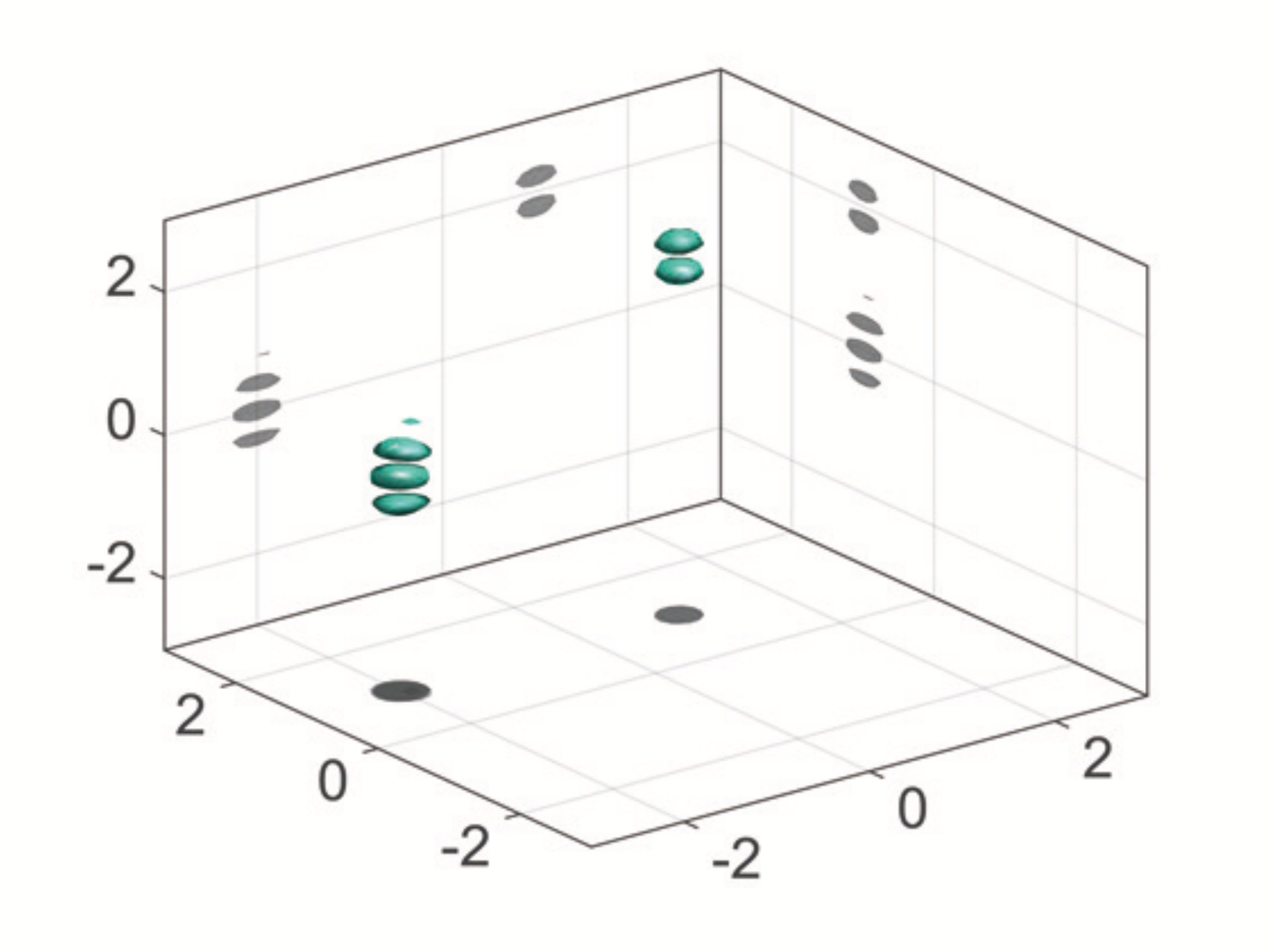}}
	\caption{Reconstruction of 3D multipoles using the DSM. The iso-surface plot of $z\mapsto |I_{3,\ell}(z)|/\max_{z\in G}|I_{3,\ell}(z)|$ with cut-off value 0.5. (a) $\ell=0$ (b) $\ell=1$ (c) $\ell=2$ (d) $\ell=3$. }
	\label{fig:example_3D_multipole_surface}
\end{figure}

\begin{figure}
	\centering
	\subfloat[]{\includegraphics[width=0.45\textwidth]{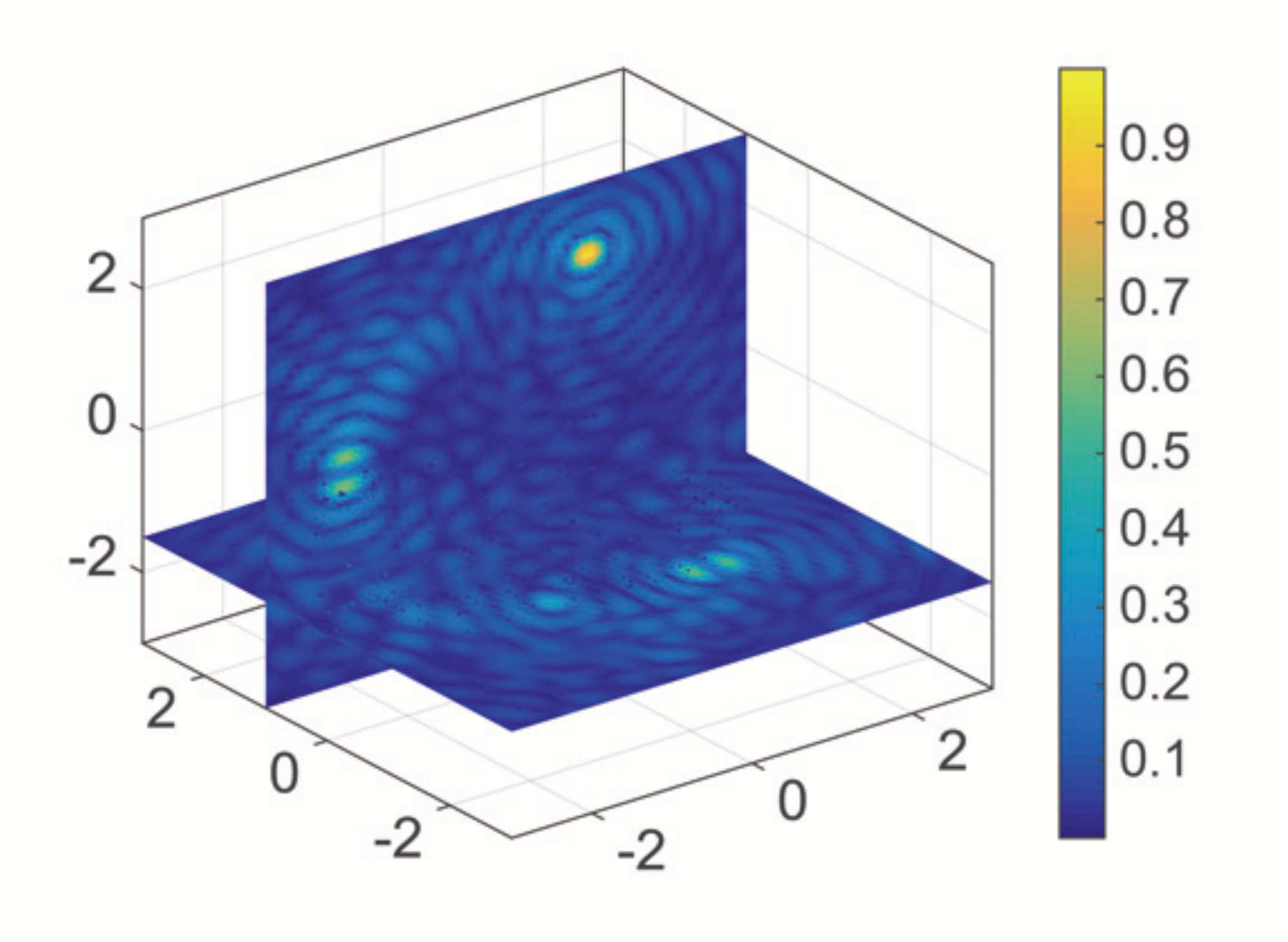}}
	\subfloat[]{\includegraphics[width=0.45\textwidth]{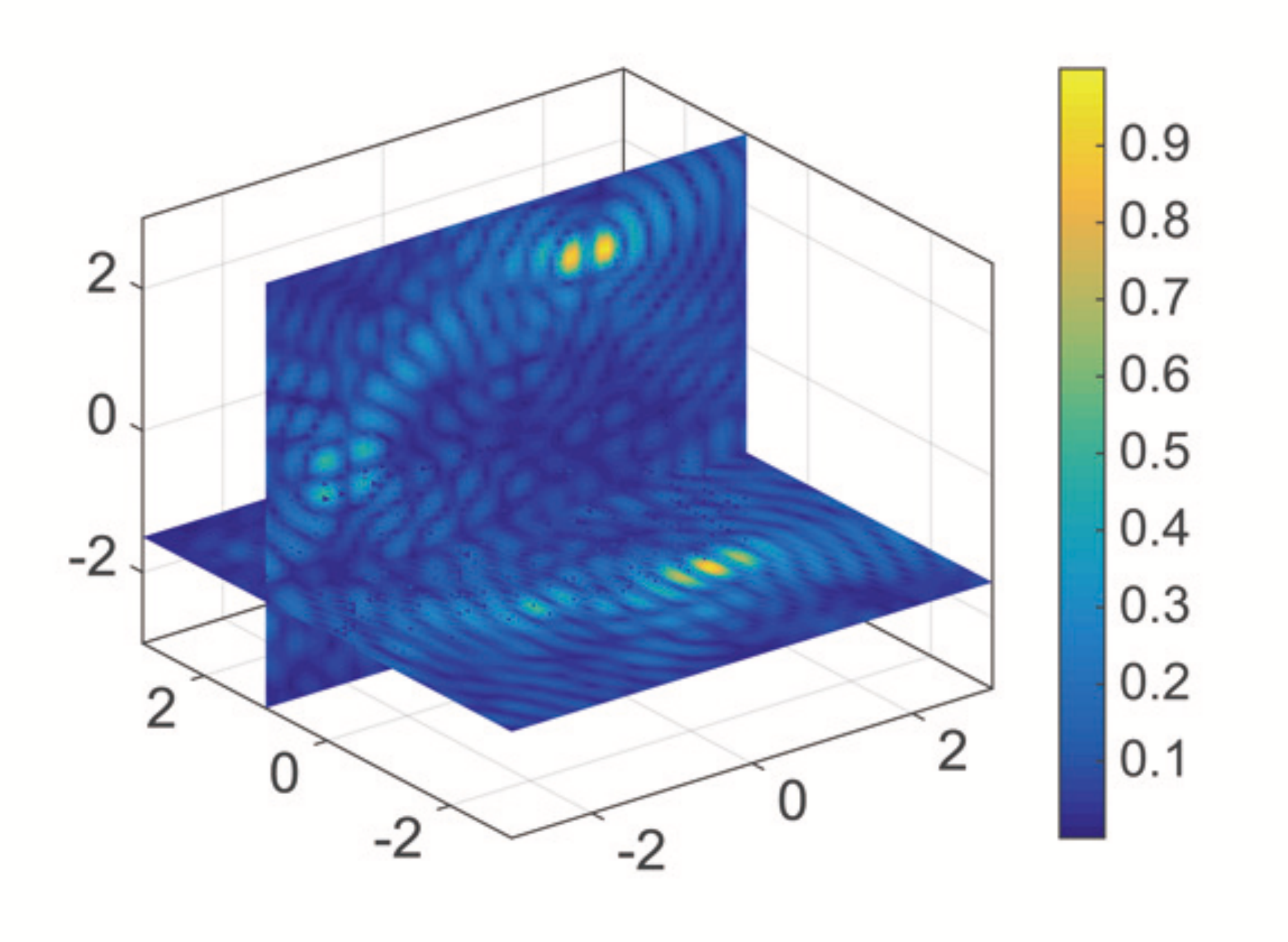}}\\
	\subfloat[]{\includegraphics[width=0.45\textwidth]{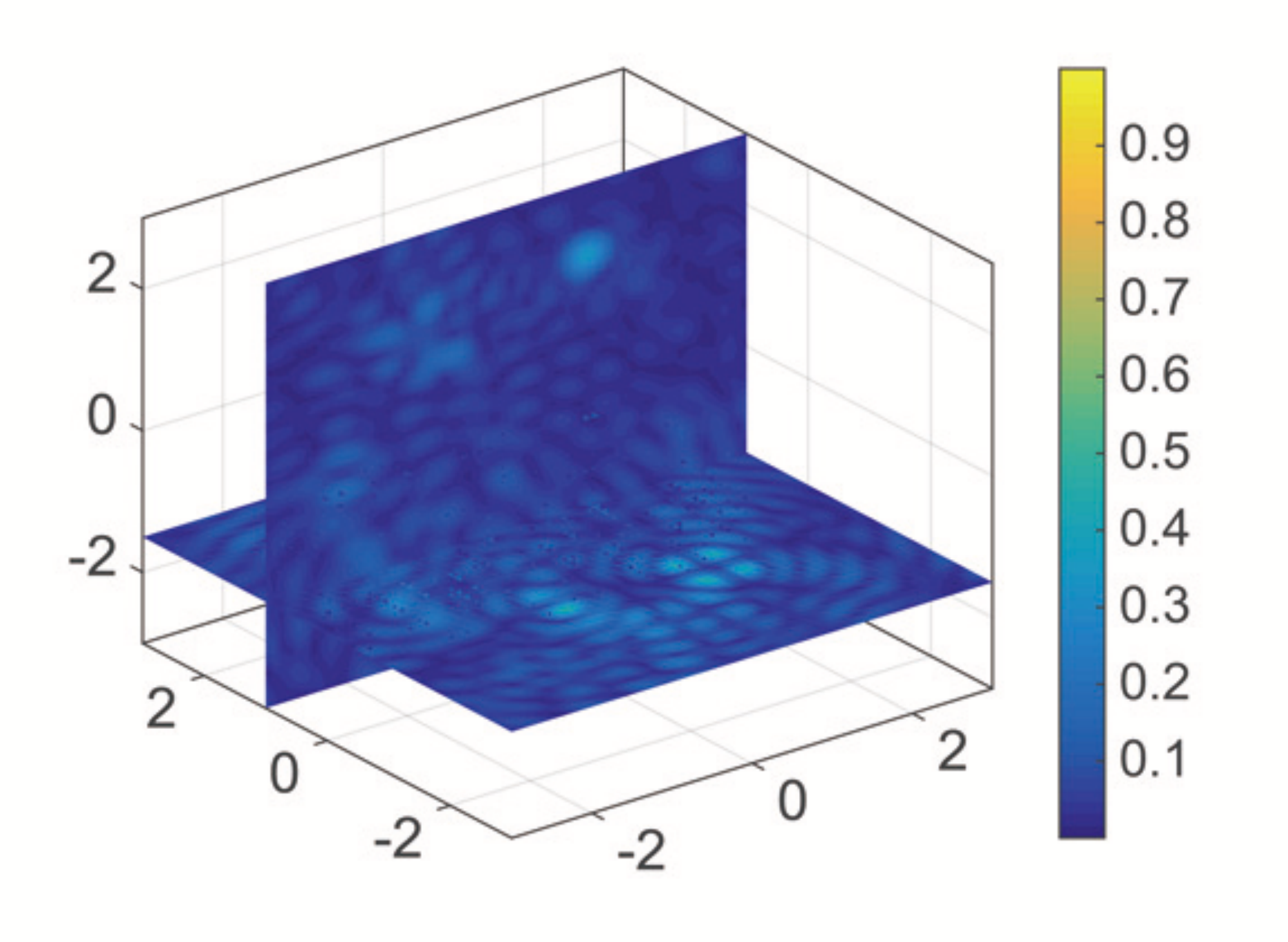}}
        \subfloat[]{\includegraphics[width=0.45\textwidth]{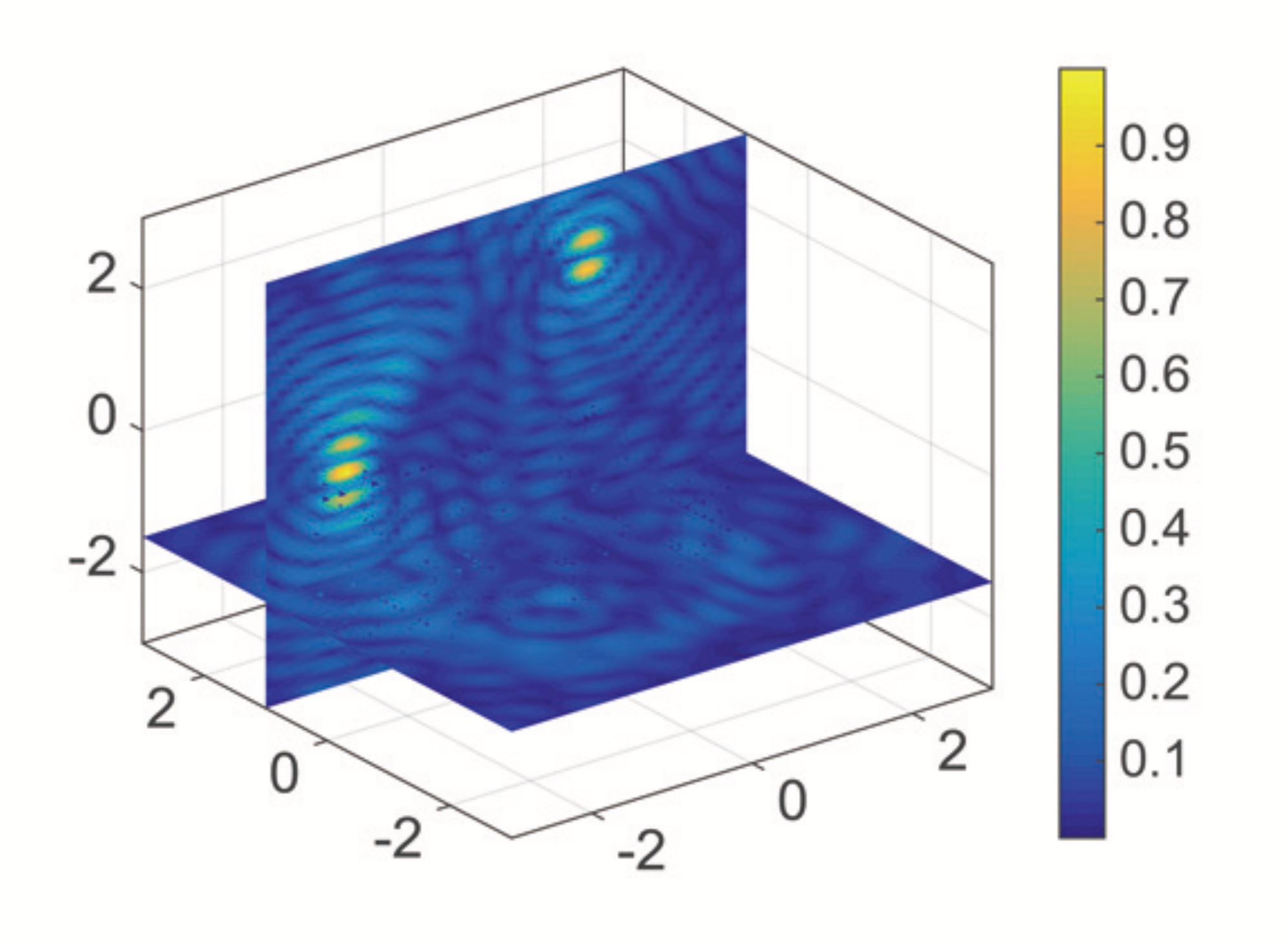}}
	\caption{Reconstruction of a monopole and two dipoles in 3D. The cross-sections of $z\mapsto |I_{3,\ell}(z)|/\max_{z\in G}|I_{3,\ell}(z)|$ are drawn with slice positions $x_2=1$ and $x_3=-1.5$. (a) $\ell=0$ (b) $\ell=1$ (c) $\ell=2$ (d) $\ell=3$.}
	\label{fig:example_3D_multipole_slice}
\end{figure}

\noindent {\bf Example 5}\quad In the last example, we will try to image the locations of a monopole and two dipoles. We just replace two of the monopoles in Example 4 by two dipoles and adjust the source intensities, so the geometry configuration of the forward problem remains the same, see Figure \ref{fig:example_3D_monopole_geometry}. The noise level is chosen as $\epsilon=15\%$ and some other parameters are listed in Table \ref{tab:example_3D_multipole}. The iso-surface and cross-section visualizations of the indicators are depicted in Figure \ref{fig:example_3D_multipole_surface} and \ref{fig:example_3D_multipole_slice}, respectively. This example illustrates that large perturbations of data (e.g. up to 15\% noise) could be well diminished by the indicator functions.

\begin{table}
\caption{\label{tab:example_3D_multipole} Reconstruction of multipoles in 3D.}
\begin{indented}\item[]
\begin{tabular}{cccc}
  \br
& & \multicolumn{2}{c}{\underline{\qquad\qquad\quad Location \qquad\qquad\qquad}} \\
Type \& Label  & Intensity & Exact & Reconstructed  \\
  \mr
monopole  & 9                 & $(1, 1, 2)$        &  $(1.0692, 1.0456, 1.9452)$ \\
dipole 1     & $(1, 0, 0)$   & $(1, -1, -1.5)$   &  $(0.8559, -0.9729, -1.5579)$ \\
dipole 2     & $(0, 0, 1)$   & $(-2, 1, 0)$       &  $(-2.0074, 0.9144, -0.0197)$ \\
  \br
\end{tabular}
\end{indented}
\end{table}

\section{Concluding remarks}

In this paper, two direct sampling schemes, DSM and DSM2, are proposed for solving the inverse source problem of locating the multipolar sources in the Helmholtz equation from the measured Cauchy data. The locating schemes are based on some novel indicator functions, which could be calculated directly via simple integrations and thus no solution process is needed. Mathematically, the indicating behaviors are rigorously analyzed. On the basis of the above theoretical justifications and the validating numerical examples, we conclude that the sampling schemes are easy to implement, robust to the measurement noise and capable of identifying the source locations effectively and efficiently.

A more general inverse source problems can be stated as: to reconstruct the source $F$ of the form \eref{FForm} in the Helmholtz equation \eref{Helmholtz} from the boundary measurements
$u|_\Gamma$ and $\partial_{\nu}u|_\Gamma$ with a single wavenumber $k$. Motivated by Theorem \ref{stability} and Remark \ref{remark_stability}, the extreme values of the indicators coincides with the source intensities as long as the exact source locations are determined. Since the reconstruction of source intensities is beyond the scope of this paper, a stable method for determining the source intensities deserves investigations in the future. We also plan to investigate the extension of the direct sampling methods to the inverse source problems for Maxwell's equations in our future works.

\appendix

\section{}

\Proof[Proof of Lemma~2.3]

Let us denote
\begin{eqnarray*}
z&=|z|\hat{z}=|z|(\sin \alpha  \cos \beta, \sin \alpha \sin \beta, \cos \alpha),\\
d &= (\sin \theta  \cos \varphi, \sin \theta \sin \varphi, \cos \theta),
\end{eqnarray*}
then from the {\it Jacobi-Anger expansion} (see \cite[p.33]{Colton})
\[
  \e^{\i k d \cdot z}=\sum_{n=0}^{\infty}\i^n(2n+1)j_n(k|z|)P_n(\cos \phi)
\]
and {\it addition theorem} (see \cite[p.27]{Colton})
\[
   \sum_{m=-n}^{n}Y_{n}^{m}(d)\overline{Y_{n}^{m}(\hat{z})}=\frac{2n+1}{4\pi}P_n(\cos \phi),
\]
we see
\begin{eqnarray} \label{3dJacobi}
  \e^{\i k d \cdot z}=4\pi\sum_{n=0}^{\infty}\sum_{m=-n}^{n}\i^n j_n(k|z|)Y_{n}^{m}(d)\overline{Y_{n}^{m}(\hat{z})},
\end{eqnarray}
where $j_n$ are the spherical Bessel functions of order $n$, $P_n$ are Legendre polynomials of degree $n$, $Y_{n}^{m}$ are the spherical harmonics, $\phi$
denotes the angle between $d$ and $z$, and the overbar denotes the complex conjugate.
By integrating \eref{3dJacobi} over $\mathbb{S}^2$ with respect to $d$ and using the definition of the spherical harmonics (see \cite[p.26]{Colton})
\begin{eqnarray} \label{harmonics}
  Y_{n}^{m}(\theta,\varphi)= \gamma_n^m P_n^{|m|}(\cos \theta)\e^{\i m \varphi},\quad  m=-n,\cdots,n,\  n=0,1,\cdots,
\end{eqnarray}
where
\[
\gamma_n^m=\sqrt{\frac{2n+1}{4\pi}\frac{(n-|m|)!}{(n+|m|)!}},
\]
and $P_n^{|m|}$ are the associated Legendre functions,  we obtain
\begin{eqnarray*}
   \int_{\mathbb{S}^2}  \e^{\i k d\cdot z}\d s(d)
   =4\pi\sum_{n=0}^{\infty}\sum_{m=-n}^{n}\i^n j_n(k|z|)\overline{Y_{n}^{m}(\hat{z})}
     \int_{\mathbb{S}^2} Y_{n}^{m}(d)\d s(d)
\end{eqnarray*}
and
\begin{eqnarray*}
  \int_{\mathbb{S}^2} Y_{n}^{m}(d)\d s(d)&=&\gamma_n^m \int_0^{2\pi}\e^{\i m \varphi} \d \varphi
      \int_0^\pi \sin \theta P_n^{|m|}(\cos \theta) \d \theta
   \\
   &=&\left\{ \begin{array}{ll}
         2\sqrt{\pi}, &\qquad m=n=0, \\
         0,  &\qquad  \mathrm{else},
       \end{array}\right.
\end{eqnarray*}
which implies
\begin{eqnarray}\label{3d0}
   \int_{\mathbb{S}^2}  \e^{\i k d\cdot z}\d s(d)=4\pi j_0(k|z|).
\end{eqnarray}

Similarly, from  the orthogonality
\[
   \int_{-1}^1   P_n^{m}(t)P_{n'}^{m}(t) \d t=\frac{2(n+m)!}{(2n+1)(n-m)!}\delta_{nn'},\quad   0\leq m \leq n, n',
\]
with the usual meaning for the Kronecker symbol $\delta_{nn'}$, we derive
\begin{eqnarray*}
  \int_{\mathbb{S}^2}d_1 Y_{n}^{m}(d)\d s(d)&=&\gamma_n^m \int_0^{2\pi}\cos \varphi \e^{\i m \varphi} \d \varphi
      \int_0^\pi \sin^2 \theta P_n^{|m|}(\cos \theta) \d \theta
   \\
     &=&\gamma_n^m \int_0^{2\pi}\cos \varphi \e^{\i m \varphi} \d \varphi
      \int_{-1}^1  \sqrt{1-t^2} P_n^{|m|}(t) \d t
   \\
     &=&\gamma_n^m \int_0^{2\pi}\cos \varphi \e^{\i m \varphi} \d \varphi
      \int_{-1}^1  P_1^1(t) P_n^{|m|}(t) \d t
   \\
   &=&\left\{ \begin{array}{ll}
          \displaystyle\frac{4\pi}{3}\gamma_1^1, &\quad m=\pm 1,\ n=1, \\
         0,  &\quad  \mathrm{ else },
       \end{array}\right.
\end{eqnarray*}
\begin{eqnarray*}
  \int_{\mathbb{S}^2}d_2 Y_{n}^{m}(d)\d s(d)
     &=&\gamma_n^m \int_0^{2\pi}\sin\varphi \e^{\i m \varphi} \d \varphi
      \int_{-1}^1  P_1^1(t) P_n^{|m|}(t) \d t
   \\
   &=&\left\{ \begin{array}{ll}
         \pm \displaystyle\frac{4\pi\i}{3}\gamma_1^1, &\quad m=\pm 1,\  n=1, \\
         0,  &\quad  \mathrm{ else },
       \end{array}\right.
\end{eqnarray*}
\begin{eqnarray*}
  \int_{\mathbb{S}^2}d_3 Y_{n}^{m}(d)\d s(d)
     &=&\gamma_n^m \int_0^{2\pi} \e^{\i m \varphi} \d \varphi
      \int_{-1}^1  P_1(t) P_n^{|m|}(t) \d t
   \\
   &=&\left\{ \begin{array}{ll}
          \displaystyle\frac{4\pi}{3}\gamma_1^0, &\quad m=0,\ n=1, \\
         0,  &\quad  \mathrm{ else },
       \end{array}\right.
\end{eqnarray*}
\begin{eqnarray*}
  \int_{\mathbb{S}^2}d_1^2 Y_{n}^{m}(d)\d s(d)&=&\gamma_n^m \int_0^{2\pi}\cos^2 \varphi \e^{\i m \varphi} \d \varphi
      \int_0^\pi \sin^3 \theta P_n^{|m|}(\cos \theta) \d \theta
   \\
     &=&\gamma_n^m \int_0^{2\pi}\frac{1+\cos2\varphi}{2} \e^{\i m \varphi} \d \varphi
      \int_{-1}^1  (1-t^2) P_n^{|m|}(t) \d t
    \\
     &=&\frac{\gamma_n^m}{6} \int_0^{2\pi}(1+\cos2\varphi) \e^{\i m \varphi} \d \varphi
      \int_{-1}^1  P_2^2(t) P_n^{|m|}(t) \d t
    \\
   &=&\left\{ \begin{array}{ll}
          \displaystyle\frac{4\pi}{3}\gamma_0^0, &\quad m=0,\ n=0, \\
         - \displaystyle\frac{4\pi}{15}\gamma_2^0, &\quad m=0,\ n=2, \\
           \displaystyle\frac{8\pi}{5}\gamma_2^2, &\quad m=\pm 2,\ n=2, \\
         0,  &\quad  \mathrm{else},
       \end{array}\right.
\end{eqnarray*}
\begin{eqnarray*}
  \int_{\mathbb{S}^2}d_2^2 Y_{n}^{m}(d)\d s(d)
     &=&\gamma_n^m \int_0^{2\pi}\frac{1-\cos2\varphi}{2} \e^{\i m \varphi} \d \varphi
      \int_{-1}^1  (1-t^2) P_n^{|m|}(t) \d t
    \\
     &=&\frac{\gamma_n^m}{6} \int_0^{2\pi}(1-\cos2\varphi) \e^{\i m \varphi} \d \varphi
      \int_{-1}^1  P_2^2(t) P_n^{|m|}(t) \d t
    \\
   &=&\left\{ \begin{array}{ll}
          \displaystyle\frac{4\pi}{3}\gamma_0^0, &\quad m=0,\ n=0, \\
         - \displaystyle\frac{4\pi}{15}\gamma_2^0, &\quad m=0,\ n=2, \vspace{1mm}\\
          - \displaystyle\frac{8\pi}{5}\gamma_2^2, &\quad m=\pm 2,\ n=2, \\
         0,  &\quad  \mathrm{else},
       \end{array}\right.
\end{eqnarray*}
\begin{eqnarray*}
  \int_{\mathbb{S}^2}d_3^2 Y_{n}^{m}(d)\d s(d)
     &=&\gamma_n^m \int_0^{2\pi} \e^{\i m \varphi} \d \varphi
      \int_{-1}^1  t^2 P_n^{|m|}(t) \d t
    \\
     &=&\left\{ \begin{array}{ll}
          \displaystyle\frac{4\pi}{3}\gamma_0^0, &\quad m=0,\ n=0, \vspace{1mm}\\
          \displaystyle\frac{8\pi}{15}\gamma_2^0, &\quad m=0,\ n=2, \\
         0,  &\quad  \mathrm{ else },
       \end{array}\right.
\end{eqnarray*}
\begin{eqnarray*}
  \int_{\mathbb{S}^2}d_1d_2 Y_{n}^{m}(d)\d s(d)
     &=&\frac{\gamma_n^m}{6} \int_0^{2\pi}\sin2\varphi \e^{\i m \varphi} \d \varphi
      \int_{-1}^1  P_2^2(t) P_n^{|m|}(t) \d t
    \\
   &=&\left\{ \begin{array}{ll}
          \pm \displaystyle\frac{8\pi\i}{5}\gamma_2^2, &\quad m=\pm 2,\ n=2, \\
         0,  &\quad  \mathrm{ else },
       \end{array}\right.
\end{eqnarray*}
\begin{eqnarray*}
  \int_{\mathbb{S}^2}d_1d_3 Y_{n}^{m}(d)\d s(d)
     &=&\frac{\gamma_n^m}{3} \int_0^{2\pi}\cos\varphi \e^{\i m \varphi} \d \varphi
      \int_{-1}^1  P_2^1(t) P_n^{|m|}(t) \d t
    \\
   &=&\left\{ \begin{array}{ll}
          \displaystyle\frac{4\pi}{5}\gamma_2^1, &\quad m=\pm 1,\ n=2, \\
         0,  &\quad  \mathrm{ else },
       \end{array}\right.
\end{eqnarray*}
and
\begin{eqnarray*}
  \int_{\mathbb{S}^2}d_2d_3 Y_{n}^{m}(d)\d s(d)
     &=&\frac{\gamma_n^m}{3} \int_0^{2\pi}\sin\varphi \e^{\i m \varphi} \d \varphi
      \int_{-1}^1  P_2^1(t) P_n^{|m|}(t) \d t
    \\
   &=&\left\{ \begin{array}{ll}
          \pm\displaystyle\frac{4\pi\i}{5}\gamma_2^1, &\quad m=\pm 1,\ n=2, \\
         0,  &\quad  \mathrm{ else },
       \end{array}\right.
\end{eqnarray*}
which yield
\begin{eqnarray}
   \quad\int_{\mathbb{S}^2} d_1 \e^{\i k d\cdot z}\d s(d)
   &=& 4\pi\i j_1(k|z|)|\sin \alpha |\cos\beta, \label{3d1} \\
   \quad\int_{\mathbb{S}^2} d_2 \e^{\i k d\cdot z}\d s(d)
   &=& 4\pi\i j_1(k|z|) |\sin \alpha |\sin\beta, \label{3d2} \\
   \quad \int_{\mathbb{S}^2} d_3 \e^{\i k d\cdot z}\d s(d)
   &=&4\pi\i  j_1(k|z|)\cos\alpha, \label{3d3} \\
   \quad\int_{\mathbb{S}^2} d_1^2 \e^{\i k d\cdot z}\d s(d)
   &=&\frac{4\pi}{3}j_0(k|z|)+\frac{2\pi}{3}j_2(k|z|)(3\cos^2\alpha-1) \nonumber \\
   &&-\frac{2\pi}{3}j_2(k|z|)\sin^2\alpha\cos2\beta, \label{3d4} \\
   \quad\int_{\mathbb{S}^2} d_2^2 \e^{\i k d\cdot z}\d s(d)
   &=&\frac{4\pi}{3}j_0(k|z|)+\frac{2\pi}{3}j_2(k|z|)(3\cos^2\alpha-1) \nonumber
   \\\label{3d5}
   &&+\frac{2\pi}{3}j_2(k|z|)\sin^2\alpha\cos2\beta ,
  \\\label{3d6}
   \quad\int_{\mathbb{S}^2} d_3^2 \e^{\i k d\cdot z}\d s(d)
   &=&\frac{4\pi}{3}j_0(k|z|)-\frac{4\pi}{3}j_2(k|z|)(3\cos^2\alpha-1),
  \\\label{3d7}
   \int_{\mathbb{S}^2} d_1d_2 \e^{\i k d\cdot z}\d s(d)
   &=&-2\pi j_2(k|z|)\sin^2\alpha\sin2\beta ,
  \\\label{3d8}
   \int_{\mathbb{S}^2} d_1d_3 \e^{\i k d\cdot z}\d s(d)
   &=&-4\pi j_2(k|z|)|\sin \alpha|\cos\alpha\cos2\beta ,
  \\\label{3d9}
   \int_{\mathbb{S}^2} d_2d_3 \e^{\i k d\cdot z}\d s(d)
   &=&-4\pi j_2(k|z|)|\sin \alpha|\cos\alpha\sin2\beta .
\end{eqnarray}

Let $h^{(1)}_{n}$ be the spherical Hankel function of the first kind of order $n$. Then from the asymptotic behavior (see \cite[p.31]{Colton})
\[
   h^{(1)}_{n}(t)= \frac{1}{t}\e^{\i  (t-\frac{1}{2}n\pi-\frac{1}{2}\pi)}
 \left\{1+\mathcal{O}(t^{-1})\right\}, \quad t\rightarrow \infty,
\]
and $j_n(t)=\mathrm{Re}\{h^{(1)}_{n}(t)\}$, we obtain the estimate
 \eref{3dAsymptotic}.  \qedhere

%%%%%%%%%%%%%%%%%%%%%%%%%%%%%%%%%%%%%%%%%%%%%%%%%%
\section*{Acknowledgements}
%%%%%%%%%%%%%%%%%%%%%%%%%%%%%%%%%%%%%%%%%%%%%%%%%%

%The authors would like to thank the referees for their careful
%reading and valuable comments which improved the quality of our
%submitted manuscript, especially on choosing the source conditions,
%determining the threshold of the discretization parameter and
%testing the cases with non-smooth boundary functions.

The work of D. Zhang was supported by the NSF grants of China under 11271159 and 11671170. The work of Y. Guo was supported by the NSF grants of China under 11601107, 11671111 and 41474102.
The work of J. Li was supported by the NSF of China under the grant No.\, 11571161 and the Shenzhen Sci-Tech Fund No. JCYJ20160530184212170.

% The authors would also like to thank the anonymous referees for many constructive comments and suggestions.

\end{document}